\documentclass[12pt]{article}

\usepackage{amssymb}
\usepackage{amsbsy}
\newcommand{\no}{\noindent}
\newcommand{\be}{\begin{equation}}
\newcommand{\ee}{\end{equation}}
\newtheorem{prop}{Proposition}[section]
\newtheorem{defn}[prop]{Definition}
\newtheorem{lem}[prop]{Lemma}
\newtheorem{cor}[prop]{Corollary}

\def\Z{\mathbb Z}
\def\L{\mathcal L}
\def\P{\mathbb P}
\def\Q{\mathbb Q}
\def\H{\mathcal H}
\def\T{\mathbb T}
\def\R{\mathbb R}
\def\C{\mathbb C}
\def\K{\mathcal K}
\def\k{\mathbb K}
\def\B{\mathcal B}
\def\F{\mathcal F}
\def\G{\mathcal G}
\def\X{\mathcal X}
\def\h{\mathbb H}
\def\p{\mathcal P}

\title{Twisted K-theory}
\author{Michael Atiyah\\ \small School of Mathematics, Edinburgh University,\\ \small James Clerk Maxwell Building,\\ \small King's Buildings, Edinburgh EH9 3JZ, U.K. \and Graeme Segal\\ \small All Souls College, Oxford OX1 4AL, U.K.}
\date{3 October 2003}

\begin{document}

\maketitle

{\bf Abstract} Twisted complex $K$-theory can be defined for a space $X$ equipped with a
bundle of complex projective spaces, or, equivalently, with a bundle of
C$^*$-algebras. Up to equivalence, the twisting corresponds to an element
of $H^3(X;\Z)$. We give a systematic account of the definition and basic
properties of the twisted theory, emphasizing some points where it
behaves differently from ordinary $K$-theory. (We omit, however, its
relations to classical cohomology, which we shall treat in a sequel.)
We develop an equivariant version of the theory for the action of a
compact Lie group, proving that then the twistings are classified by
the equivariant cohomology group $H^3_G(X;\Z)$. We also consider some
basic examples of twisted $K$-theory classes, related to those appearing
in the recent work of Freed-Hopkins-Teleman.

\section{Introduction}

In classical cohomology theory the best known place where one encounters twisted coefficients is the Poincar\'{e} duality theorem, which, for a compact oriented $n$-dimensional manifold $X$, relates to the pairing between cohomology classes in complementary dimensions given by multiplication followed by integration over $X$:
	\[H^p(X;\Z)\times H^{n-p}(X:\Z)\rightarrow H^n(X;\Z)\rightarrow \Z.\]
If $X$ is not orientable there is a local coefficient system $\omega$ on $X$ whose fibre $\omega_x$ at each point $x$ is non-canonically $\Z$, and the duality pairing is
	\[H^p(X;\Z)\times H^{n-p}(X;\omega)\rightarrow H^n(X;\omega)\rightarrow\Z.\]
The difference between elements of $H^n(X;\Z)$ and of $H^n(X;\omega)$ is the difference between $n$-forms and densities.

In $K$-theory the Poincar\'{e} pairing involves twisting even when $X$ is oriented.  Let us, for simplicity, take $X$ even dimensional and Riemannian.  Then the analogue of the local system $\omega$ is the bundle $C$ of finite dimensional algebras on $X$ whose fibre $C_x$ at $x$ is the complex Clifford algebra of the cotangent space $T^*_x$ at $x$.  Alongside the usual $K$-group $K^0(X)$ formed from the complex vector bundles on $X$ there is the group $K^0_C(X)$ formed from $C$-modules, i.e. finite dimensional complex vector bundles $E$ on $X$ such that each fibre $E_x$ has an action of the algebra $C_x$.  On the sections of such a $C$-module $E$ there is a Dirac operator
	\[D_E=\Sigma\gamma_iD_i\]
(defined by choosing a connection in $E$; here $D_i$ is covariant differentiation in the $i^{\rm th}$ coordinate direction, and $\gamma_i$ is Clifford multiplication by the dual covector). In fact the module $E$ automatically has a decomposition $E=E^+ \oplus E^-$, and $D_E$ maps the space of sections $\Gamma (E^+ )$ to $\Gamma (E^- )$, and vice versa. Each component
\[D_E^{\pm}:\Gamma (E^{\pm})\rightarrow \Gamma(E^{\mp})\]
is Fredholm, and associating to $E$ the index of $D_E^+$ defines a homomorphism 
	\[K^0_C(X)\rightarrow\Z\]
which is the $K$-theory analogue of the integration map
	\[H^n(X;\omega)\rightarrow\Z.\]
Tensoring $C$-modules with ordinary vector bundles then defines the Poincar\'{e} pairing
\renewcommand{\theequation}{1.\arabic{equation}}
	\be K^0(X)\times K^0_C(X)\rightarrow K^0_C(X)\rightarrow\Z.\ee

We can define a twisted $K$-group $K^0_A$ for any bundle $A$ of finite dimensional algebras on $X$.  The interesting case is when each fibre $A_x$ is a full complex matrix algebra: equivalence classes of such bundles $A$ correspond, as we shall see, to the torsion elements in $H^3(X;\Z)$.  The class of the bundle $C$ of Clifford algebras of an even-dimensional orientable real vector bundle $E$ is the integral third Stiefel-Whitney class $W_3(E)\in H^3(X;\Z)$, the image of $w_2(E)\in H^2(X;\Z/2)$ by the Bockstein homomorphism.  In this paper we shall consider a somewhat more general class of twistings parametrized by elements of $H^3(X;\Z)$ which need not be of finite order.  From one viewpoint the new twistings correspond to bundles of infinite dimensional algebras on $X$.

In fact the bundle $C$ of Clifford algebras on a manifold $X$ is a mod 2 graded algebra, and the definition of $K^0_C$ should take the grading into account.  When this is done the pairing (1.1) expresses Poincar\'{e} duality even when $X$ is not orientable.
\bigskip

The existence of the twisted $K$-groups has been well-known to experts since the early days of $K$-theory (cf. Donovan-Karoubi [DK], Rosenberg [R]), but, having until recently no apparent role in geometry, they attracted little attention.  The rise of string theory has changed this.  In string theory space-time is modelled by a new kind of mathematical structure whose ``classical limit" is not just a Riemannian manifold, but rather one equipped with a so-called {\it B-field} [S4].  A $B$-field $\beta$ on a manifold $X$ is precisely what is needed to define a twisted $K$-group $K^0_\beta (X)$, and the elements of this group represent geometric features of the stringy space-time.  If the field $\beta$ is realized by a bundle $A$ of algebras on $X$ then $K^0_\beta (X)$ is the $K$-theory of the non-commutative algebra of sections of $A$, and it is reasonable to think of the stringy space-time as the ``non-commutative space"  --- in the sense of Connes [C] --- defined by this algebra. Many papers have appeared recently discussing twisted $K$-theory in relation to string theory, the most comprehensive account probably being that of the Adelaide school [BCMMS]. Another account is given in [TXL]. We refer to [Mo] for a physicist's approach.
\bigskip

A purely mathematical reason for being interested in twisted $K$-theory is the beautiful theorem proved recently by Freed, Hopkins, and Teleman which expresses the Verlinde ring of projective representations of the loop group $\L G$ of a compact Lie group $G$ --- a ring under the subtle operation of ``fusion" --- as a twisted equivariant $K$-group of the space $G$.  Here the twisting corresponds to the ``level", or projective cocycle, of the representations being considered.
\bigskip

In this paper we shall set out the basic facts about twisted $K$-theory simply but carefully.  There are at least two ways of defining the groups, one in terms of families of Fredholm operators, and the other as the algebraic $K$-theory of a non-commutative algebra.  We shall adopt the former, but shall sketch the latter too.  The equivariant version of the theory is of considerable interest, but it has seemed clearest to present the non-equivariant theory first, using arguments designed to generalize, and only afterwards to explain the special features of the equivariant case.
\bigskip

            The plan of the paper is as follows.
\bigskip

Section 2 discusses the main properties of bundles of infinite dimensional projective spaces, which are the ``local systems" which we shall use to define twisted K-theory.
\bigskip

Section 3 gives the definition of the twisted $K$-theory of a space $X$ equipped with a bundle $P$ of projective spaces, first as the group of homotopy classes of sections of a bundle on $X$ whose fibre at $x$ is the space of Fredholm operators in the fibre $P_x$ of $P$, then as the algebraic $K$-theory of a $C^*$-algebra associated to $X$ and $P$. The twistings by bundles of projective spaces are not the most general ones suggested by algebraic topology, and at the end of this section we mention the general case.
\bigskip

Section 4 outlines the algebraic-topological properties of twisted K-theory. The relation of the twisted theory to classical cohomology will be discussed in a sequel to this paper.
\bigskip

Section 5 describes some interesting examples of projective bundles and families of Fredholm operators in them, related to the ones occurring in the work of Freed, Hopkins, and Teleman [FHT]. In fact these are naturally equivariant examples. They have also been discussed by Mickelsson [M] (cf. also [CM]).
\bigskip

Section 6 turns to the equivariant theory, explaining the parts which are not just routine extensions of the non-equivariant discussion.
\bigskip

           Apart from that there are three technical appendices concerned with points of functional analysis with which we did not want to hold up the main text. The third is an equivariant version of Kuiper's proof of the contractibility of the general linear group of Hilbert space with the norm topology.
\bigskip

In a subsequent paper we shall discuss the relation of twisted $K$-theory to cohomology. We shall examine the effect of twisting on the Atiyah-Hirzebruch spectral sequence, on the Chern classes, and on the Chern character. We shall also see how twisting interacts with the operations in K-theory, such as the exterior powers and the Adams operations.

\section{Bundles of projective spaces}
\renewcommand{\theequation}{2.\arabic{equation}}\renewcommand{\thefootnote}{*}

The ``local systems" which we shall use to define twisted $K$-theory are bundles of infinite dimensional complex projective spaces.  This section treats their basic properties.

We shall consider locally trivial bundles $P\rightarrow X$ whose fibres $P_x$ are of the form $\P(\H)$, the projective space of a separable complex Hilbert space $\H$ which will usually, but not invariably, be infinite dimensional (we shall at least require that it has dimension $\geq 1$, so that $\P(\H)$ is non-empty).  We shall assume that our base-spaces $X$ are metrizable, though this could easily be avoided by working in the category of compactly generated spaces.  The projective-Hilbert structure of the fibres is supposed to be given.  This means that $P$ is a fibre bundle whose structural group is the projective unitary group $PU(\H)$ with the compact-open topology.\footnote{An account of the compact-open topology can be found in Appendix 1.}  The significance of this topology is that a map $X\rightarrow PU(\H)$ is the same thing as a bundle isomorphism 
	\[X\times\P(\H)\rightarrow X\times\P(\H).\]  
In fact, essentially by the Banach-Steinhaus theorem, the same is true if $PU(\H)$ has the slightly coarser topology of pointwise convergence, which is called the ``strong operator topology" by functional analysts.

Let us stress that we do not always want to assume that the structural group of our bundles is $PU (\H)$ with the {\it norm} topology, i.e. that there is a preferred class of local trivializations between which the transition functions are norm-continuous, for doing so would exclude most naturally arising bundles.  For example, if $Y\rightarrow X$ is a smooth fibre bundle with compact fibres $Y_x$ then the Hilbert space bundle $E$ on $X$ whose fibre $E_x$ is the space of $L^2$ half-densities on $Y_x$ does not admit $U(\H)$ with the norm topology as structure group, for the same reason that if $\H=L^2(G)$ is the regular representation of a group $G$ the action map $G\rightarrow U(\H)$ is not norm-continuous, even if $G$ is compact. Nevertheless, it follows from Proposition 2.1(ii) below that for many purposes we lose nothing by working with norm-continuous projective bundles, and it is simpler to do so.
\bigskip

When we have a bundle $P\rightarrow X$ of projective spaces we can construct another bundle End$(P)$ on $X$ whose fibre at $x$ is the vector space End$(\H_x)$ of endomorphisms of a Hilbert space $\H_x$ such that $P_x={\mathbb P}(\H_x)$.  For, although $\H_x$ is not determined canonically by the projective space $P_x$, if we make another choice $\tilde{\H}_x$ with ${\mathbb P}(\tilde{\H}_x)={\mathbb P}(\H_x)$ then End$(\tilde{\H}_x)$ is {\it canonically} isomorphic to End$(\H_x)$, and it makes sense to define
	\[{\rm End}(P_x)=\ {\rm End}(\H_x)=\ {\rm End}(\tilde{\H}_x).\]
This observation will play a basic role for us, and we shall use several variants of it, replacing End$(\H_x)$ by, for example, the subspaces of compact, Fredholm, Hilbert-Schmidt, or unitary operators in End$(\H_x)$.  We must beware, however, that if the structural group of $P$ does not have the norm topology we must use the compact-open topology on the fibres of End$(P)$, Fred$(P)$, or $U(P)$.  In the case of the compact or Hilbert-Schmidt operators there is no problem of this kind, for, as is proved in Appendix 1, the group $U(\H)$ with the compact-open topology acts continuously on the Banach space $\K(\H)$ of compact operators and the Hilbert space $\H^*\otimes\H$ of Hilbert-Schmidt operators.
\bigskip

Each bundle $P\rightarrow X$ of projective spaces has a class $\eta_P\in H^3(X;\Z)$ defined as follows.  Locally $P$ arises from a bundle of Hilbert spaces on $X$, so we can choose an open covering $\{X_\alpha\}$ of $X$ and isomorphisms $P\vert X_\alpha\cong\P(E_\alpha)$, where $E_\alpha$ is a Hilbert space bundle on $X_\alpha$.  If the covering $\{X_\alpha\}$ is chosen sufficiently fine \footnote{This is a slight oversimplification. Most spaces of interest posses arbitrarily fine open coverings $\{X_\alpha\}$ such that the intersections $X_{\alpha\beta}$ are contractible, and then the maps $g_{\alpha\beta}$ can be lifted to vector bundle isomorphisms, e.g. by fixing the phase of some matrix element (which is continuous in the compact-open topology). But in general we must use the standard technology of sheaf theory, which takes a limit over coverings rather than using a particular covering.} the transition functions between these ``charts" can be realized by isomorphims
	\[g_{\alpha\beta}:E_\alpha\vert X_{\alpha\beta}\rightarrow E_\beta\vert X_{\alpha\beta},\]
where $X_{\alpha\beta}=X_\alpha\cap X_\beta$, which are projectively coherent, so that over each triple intersection $X_{\alpha\beta\gamma}=X_\alpha\cap X_\beta\cap X_\gamma$ the composite
	\[g_{\gamma\alpha}g_{\beta\gamma}g_{\alpha\beta}\]
is multiplication by a circle-valued function $f_{\alpha\beta\gamma}:X_{\alpha\beta\gamma}\rightarrow\T$.  These functions $\{f_{\alpha\beta\gamma}\}$ constitute a cocycle defining an element $\tilde{\eta}_P$ of the \v{C}ech cohomology group $H^2(X;sh(\T))$, where $sh(\T)$ denotes the sheaf of continuous $\T$-valued functions on $X$.  Using the exact sequence
	\[0\rightarrow sh(\Z)\rightarrow sh(\R)\rightarrow sh(\T)\rightarrow 0\]
we can define $\eta_P$ as the image of $\tilde{\eta}_P$ under the coboundary homomorphism
	\[H^2(X;sh(\T))\rightarrow H^3(X;sh(\Z))=H^3(X;\Z)\]
(which is an isomorphism because $H^i(X;sh(\R))=0$ for $i>0$ by the existence of partitions of unity).
\bigskip

Before stating the main result of this section let us notice that bundles of projective Hilbert spaces can be tensored: the fibre $(P_1\otimes P_2)_x$ is the Hilbert space tensor product $P_{1,x}\otimes P_{2,x}$, i.e. the projective space of the Hilbert space of Hilbert-Schmidt operators $E^*_{1,x}\rightarrow E_{2,x}$, where $E^*_{i,x}$ is the dual space of $E_{i,x}$, and $\P (E_{i,x})\cong P_{i,x}$.   Furthermore, for any bundle $P$ there is a dual projective bundle $P^*$ whose points are the closed hyperplanes in $P$, and $P^*\otimes P$ comes from a vector bundle.  In fact $P^*\otimes P=\P(E)$, where $E$ is the bundle of Hilbert-Schmidt endomorphisms of $P$.  (This is a first application of the observation above that the vector space End$(\H)$ is functorially associated to the projective space $\P(\H)$, even though $\H$ itself is not.)

\begin{prop}  
\begin{enumerate}
\item[(i)] We have $\eta_P=0$ if and only if the bundle $P$ of projective spaces comes from a vector bundle $E$ on $X$.
\item[(ii)] Each element of $H^3(X;\Z)$ arises from a bundle $P$, even from one whose structure group is $PU(\H)$ with the norm topology.
\item[(iii)] If the fibres of $P$ are infinite dimensional and separable then  $P$ is determined up to isomorphism by $\eta_P$.
\item[(iv)] If $P$ has finite dimensional fibres $\P(\C^n)$ then $n\eta_P=0$.
\item[(v)] Every torsion element of $H^3 (X;\Z)$ arises from a finite dimensional bundle $P$, though a class of order $n$ need not arise from a bundle with fibre $\P(\C^n)$.
\item[(vi)] If $P_0\rightarrow P$ is a tame embedding of projective bundles, in the sense explained below, then $\eta_{P_0}=\eta_P$.  In particular, if $P$ has a continuous section then $\eta_P=0$, and if $\P$ is a fixed projective space then $\eta_P=\eta_{P\otimes\P}$.
\item[(vii)] We have $\eta_{P_1\otimes P_2}=\eta_{P_1}+\eta_{P_2}$.
\item[(viii)] We have $\eta_{P^*}=- \eta_P$.
\end{enumerate}
\end{prop}

In (vi) above, a {\it tame} embedding means one which is locally isomorphic (on $X$) to the inclusion of $X\times\P(\H_0)$ in $X\times\P(\H)$, where $\H_0$ is a closed subspace of $\H$.  A typical example of a non-tame embedding is the following.  Let $\H$ be the standard Hilbert space $L^2(0,1)$.  Then in the trivial bundle $X\times\H$ on the closed interval $X=[0,\frac{1}{2}]$ the subbundle whose fibre at $x$ is $L^2(x,1)$ is not tame.
\bigskip

Proposition (2.1), whose proof is given below, tells us that the group of isomorphism classes of projective bundles (with infinite dimensional separable fibres) under the tensor product is precisely $H^3(X;\Z)$.  We also need to know about the automorphism groups of these bundles.  An automorphism $\alpha :P\rightarrow P$ defines a complex line bundle $L_\alpha$ on $X$: the non-zero elements of the fibre of $L_\alpha$ at $x$ are the linear isomorphisms $E_x\rightarrow E_x$ which induce $\alpha\vert P_x$, where $P_x=\P(E_x)$.  (We have already pointed out that the choice of $E_x$ is irrelevant.)

\begin{prop}
For a projective bundle $P$ with infinite dimensional separable fibres the assignment $\alpha\mapsto L_\alpha$ identifies the group of connected components  of the automorphism group of $P$ with the group $H^2(X;\Z)$ of isomorphism classes of complex line bundles on $X$.
\end{prop}

The proof will be given presently.
\bigskip

\noindent {\bf Proof of Proposition (2.1)}
\bigskip

(i) This is immediate because the vanishing of the \v{C}ech cohomology class $\tilde{\eta}_P\in H^2(X;sh(\T))$ defined by transition functions $\{g_{\alpha\beta}\}$ is precisely the condition that the $g_{\alpha\beta}$ can be multiplied by functions $\lambda_{\alpha\beta}:X_{\alpha\beta}\rightarrow \T$ to make them exactly coherent. 
\bigskip

(ii) Because the unitary group $U(\H)$ of an infinite dimensional Hilbert space is contractible --- with either the norm topology, or the compact-open topology (see Appendix 2) --- the projective group $PU(\H)$ has the homotopy type of an Eilenberg-Maclane space $K(\Z,2)$, and its classifying space $BPU(\H)$ is accordingly a $K(\Z,3)$.  Thus any element of $H^3(X;\Z)$ corresponds to a map $f:X\rightarrow BPU(\H)$, and hence to the bundle of projective spaces pulled back by $f$ from the universal bundle on $BPU(\H)$.
\bigskip

(iii) Any bundle can be pulled back from the universal bundle, and homotopic maps pull back isomorphic bundles.
\bigskip

(iv) The commutative diagram of exact sequences
	\[\begin{array}{ccccc} \mu_n & \longrightarrow & SU_n & \longrightarrow & PU_n\\ \downarrow && \downarrow && \downarrow\\ \T & \longrightarrow & U_n & \longrightarrow & PU_n,\end{array}\]
where $\mu_n$ is the group of $n^{\rm th}$ roots of unity, and the right-hand vertical map is the identity, shows that the invariant $\tilde{\eta}_P\in H^2(X;sh(\T))$, when $P$ has structural group $PU_n$, comes from $H^2(X;sh(\mu_n))$, and hence has order dividing $n$.
\bigskip

(v)	(The following argument is due to Serre, see [G].)  If $l$ divides $m$ --- say $m=lr$ --- we have an inclusion $PU_l\rightarrow PU_m$ given by tensoring with $\C^r$.  By Bott periodicity the homotopy groups $\pi_i(BPU_l)$ for $i<2l-1$ are given by
	\begin{eqnarray*} \pi_2(BPU_l)&=&\Z/l\\ \pi_i(BPU_l)&=&\Z \ \mbox{for $i$ even and $>1$}\\ \pi_i(BPU_l)&=&0 \ \mbox{for $i$ odd.}\end{eqnarray*}
The inclusion $PU_l\rightarrow PU_m$ induces multiplication by $r=m/l$ on all homotopy groups, so we have

	\begin{eqnarray*} \pi_2(BPU_\infty)&=&\Q/\Z\\ \pi_i(BPU_\infty)&=&\Q \ \mbox{for $i$ even and $>1$}\\ \pi_i(BPU_\infty)&=&0 \ \mbox{for $i$ odd.}\end{eqnarray*}
Thus $BPU_\infty$ can be constructed from the Eilenberg-Maclane space $K(\Q/\Z,2)$ by successively forming fibrations over it with fibres $K(\Q,2j)$.  A fibration with fibre $K(\Q,2j)$ on a base-space $Y$ is determined by an element of $H^{2j+1}(Y;\Q)$.  Now $K(\Q/\Z,2)$ has the rational cohomology of a point, while the other Eilenberg-Maclane spaces involved have rational cohomology only in even dimensions.  So
	\[BPU_\infty\simeq K(\Q/\Z,2)\times K(\Q,4)\times K(\Q,6)\times \ldots .\]
This means that every element $\eta '$ of $H^2(X;\Q/\Z)$ can be realized by a $BPU_\infty$-bundle $P$ whose invariant $\eta_P$ is the image of $\eta '$ in $H^3(X;\Z)$.  But from the Bockstein sequence for
	\[0\rightarrow \Z\rightarrow \Q\rightarrow \Q/\Z\rightarrow 0\]
the torsion elements of $H^3(X;\Z)$ are precisely the image of $H^2(X;\Q/\Z)$.

There is, however, no reason to expect that when $n\eta =0$ we can represent the class of $\eta$ by a bundle with fibre $\P(\C^n)$.  We have seen, for example, that the class of the bundle Cliff$(E)$ of Clifford algebras of a $2k$-dimensional real vector bundle $E$ --- or, equivalently, of the projective bundle of spinors of $E$ --- is $W_3(E)\in H^3(X;\Z)$, which is of order 2, while the projective bundle of spinors has dimension $2^k-1$, and its class need not be represented by a bundle of lower dimension.  To have a concrete counterexample we can reason as follows.  The invariant of a bundle with fibre $\P(\C^2)$ is given by a map
	\[BPU_2\rightarrow K(\Z/2,2)\rightarrow K(\Z,3).\]
If every invariant of order 2 came from a $PU_2$-bundle then the map $K(\Z/2,2)\rightarrow K(\Z,3)$ would factorize
	\[K(\Z/2,2)\rightarrow BPU_2\rightarrow K(\Z,3),\]
and taking loops would give
	\[K(\Z/2,1)\rightarrow PU_2\rightarrow K(\Z,2),\]
which is impossible because the Bockstein map $K(\Z/2,1)\rightarrow K(\Z,2)$ (i.e. $\R P^\infty\hookrightarrow\C P^\infty)$ clearly does not factorize through a finite dimensional space.
\bigskip

(vi) This follows by the argument of case (iv) from the diagram
	\[\begin{array}{lcccc} \T & \longrightarrow & U(\H_0) & \longrightarrow & PU(\H_0)\\ \uparrow\cong && \uparrow && \uparrow\\ \T & \longrightarrow & U(\H,\H_0) & \longrightarrow & PU(\H,\H_0)\\
\downarrow\cong && \downarrow &&\downarrow\\
\T & \longrightarrow & U(\H) & \longrightarrow & PU(\H)\end{array}\]
where $U(\H,\H_0)=\{u\in U(\H):u(\H_0)=\H_0\}$.
\bigskip

(vii) Here we consider
	\[\begin{array}{ccccc} \T\times\T & \longrightarrow & U(\H_1)\times U(\H_2) & \longrightarrow & PU(\H_1)\times PU(\H)\\ \downarrow && \downarrow && \downarrow\\ \T & \longrightarrow & U(\H_1\otimes\H_2) & \longrightarrow & PU(\H_1\otimes\H_2)\end{array}\]
where the left-hand vertical map is composition in $\T$.
\bigskip

(viii) This follows from (vii).
\bigskip

\noindent {\bf Proof of Proposition (2.2)}

An automorphism of $P$ is a section of a bundle on $X$ whose fibre is $PU(\H)$.  This bundle, however, comes from one with fibre $U(\H)$, and so it is trivial.  The group of automorphisms can therefore be identified with the maps from $X$ to $PU(\H)$, which is an Eilenberg-Maclane space $K(\Z,2)$.

\bigskip
\subsubsection*{Remark}

In fact the natural objects that can be used to twist $K$-theory are not simply bundles $P\rightarrow X$ of projective spaces, but rather are bundles of projective spaces in which a unitary involution is given in each fibre $P_x$.  An involution in a projective space $\P$ expresses it as the join of two disjoint closed projective subspaces $\P^+$ and $\P^-$ which, despite the notation, are {\it not} ordered.  We shall always assume that $P^+_x$ and $P^-_x$ fit together locally to form tame subbundles of $P$.  Thus the involution defines a double covering of $X$, and hence a class $\xi_P\in H^1(X;\Z/2)$.

Let Proj$^{\pm}(X)$ denote the group of isomorphism classes of infinite dimensional projective Hilbert space bundles with involution on $X$, under the operation of graded tensor product.

\begin{prop} As sets we have
	\[{\rm Proj}^{(\pm)}(X)\cong H^1(X;\Z/2)\times H^3(X;\Z)\]
canonically, but the tensor product of bundles induces the product
	\[(\xi_1,\eta_1).(\xi_2,\eta_2)=(\xi_1+\xi_2,\eta_1+\eta_2+\beta (\xi_1\xi_2))\]
on the cohomology classes, where $\xi_1\xi_2 \in H^2 (X;\Z/2)$ is the cup-product, and
	\[\beta:H^2(X;\Z/2)\rightarrow H^3(X;\Z)\]
is the Bockstein homomorphism.
\end{prop}

In other words, we have an exact sequence
	\[0\rightarrow H^3(X;\Z)\rightarrow\,{\rm Proj}^{(\pm)}(X)\rightarrow H^1(X;\Z/2)\rightarrow 0\]
which is split (because every element of ${\rm Proj}^{(\pm)}$ has order two), but not canonically split.  The Bockstein cocycle describing the extension expresses the extent to which the forgetful functor from projective spaces with involution to projective spaces does not respect the tensor product.  The proof of (2.3) is very simple.  We can think of elements of $H^1(X;\Z/2)$ as real line bundles on $X$, and can define a map
	\[H^1(X;\Z/2)\rightarrow\,{\rm Proj}^{(\pm)}(X)\]
by taking a line bundle $L$ to $\P(S_L\otimes\H)$, where $S_L$ is an irreducible graded module for the bundle of Clifford algebras $C(L)$, and $\H$ is a fixed Hilbert space.  Now
	\begin{eqnarray*} \P(S_{L_1}\otimes\H)\otimes\P(S_{L_2}\otimes\H)&\cong&\P(S_{L_1}\otimes S_{L_2}\otimes\H\otimes\H)\\ &\cong&\P(S_{L_1\oplus L_2}\otimes\H),\end{eqnarray*}
where everything is understood in the graded sense.  But
	\begin{eqnarray*} W_3(L_1\oplus L_2)&=&\beta (w_2(L_1\oplus L_2))\\ &=&\beta (w_1(L_1)w_1(L_2)),\end{eqnarray*}
which is the assertion of (2.3).

\bigskip
For simplicity, in the rest of this paper we shall not pursue this generalization, but for the most part will keep to the twistings corresponding to elements of $H^3$. The other extreme, when the twisting is given by an element of $H^1(X;\Z /2)$ alone, is a special case of the version of $K$-theory developed by Atiyah and Hopkins [AH].
\bigskip

\section{The definition}

It is well known (see [A] Appendix) that the space Fred$(\H)$ of Fredholm operators in an infinite dimensional Hilbert space $\H$, with the norm topology, is a representing space for $K$-theory, i.e. that
	\[K^0(X)\cong [X;{\rm Fred}(\H)]\]
for any space $X$, where [ \ \ ; \ \ ] denotes the set of homotopy classes of continuous maps.  The basic observation for twisting $K$-theory is that when $P$ is a bundle on $X$ with fibre $\P(\H)$ there is an associated bundle Fred$(P)$ with fibre Fred$(\H)$, and we can define $K^0_P(X)$ as the set of homotopy classes of sections of Fred$(P)$.

If the bundle $P$ admits the projective unitary group $PU(\H)_{\rm norm}$ with the {\it norm} topology as its structure group this is straightforward, as $PU(\H)_{\rm norm}$ acts on Fred$(\H)$ by conjugation.  But, as we have explained, we want to avoid that assumption. (To be quite clear, for any given projective bundle $P$ we could, by 2.1(ii), choose a reduction of the structure group to $PU(\H)_{\rm norm}$, but we could not then expect a natural family of Fredholm operators in $P$ to define a continuous section of Fred$(P)$.)  We can, of course, in any case construct a bundle whose fibre is Fred$(\H)_{\rm c.o.}$ with the compact-open topology, but Fred$(\H)_{\rm c.o.}$ does not represent $K$-theory: it is a contractible space (see Appendix 2), and the index is not a continuous function on it.

We can deal with this problem in  various ways.  The simplest is to replace Fred$(\H)$ by another representing space for $K$-theory on which $PU(\H)_{\rm c.o.}$ does act continuously.  One such space is the restricted Grassmannian Gr$_{\rm res}(\H)$ described in Chap. 7 of [PS].  In practical applications of the theory, however, $K$-theory elements are more commonly represented by families of Fredholm operators --- often elliptic differential operators --- than by maps into Grassmannians. We therefore stay with Fredholm operators, and we can do this by defining a modified space of operators, bearing in mind that a continuously varying Fredholm operator usually has a natural continuously varying parametrix.  An operator $A:\H\rightarrow\H$ is Fredholm if and only if it is invertible modulo compact operators, i.e. if there exists a ``parametrix"  $B:\H\rightarrow\H$ such that $AB-1$ and $BA-1$ are compact.  Let us therefore consider the set Fred$'(\H)$ of pairs $(A,B)$ of Fredholm operators related in this way.  Ignoring topology for the moment, notice that the projection $(A,B)\mapsto A$ makes Fred$'(\H)$ a bundle of affine spaces over Fred$(\H)$ whose fibres are isomorphic to the vector space $\K$ of compact operators.  We shall give Fred$'(\H)$ the topology induced by the embedding
	\[(A,B)\mapsto(A,B,AB-1,BA-1)\]
in $\B\times\B\times\K\times\K$, where $\B$ is the bounded operators in $\H$ with the compact-open topology and $\K$ is the compact operators with the norm topology.

A proof of the following proposition is implicit in [S3], where a more general situation is treated. But for clarity we have included a direct proof of 3.1(i) in Appendix 2, while 3.1(ii) is proved in Appendix 1.

\begin{prop}
\begin{enumerate}\item[(i)] Fred$'(\H)$ is a representing space for $K$-theory.
\item[(ii)] The group PGL$(\H)$ with the compact-open topology acts continuously on Fred$'(\H)$ by conjugation.
\end{enumerate}
\end{prop}

If $P\rightarrow X$ is an infinite dimensional bundle of projective spaces Proposition (3.1) allows us to define the associated bundle Fred$'(P)$, and we can define $K^0_P(X)$ as the group of homotopy classes of its actions.  To deal with the multiplicativity properties of $K$-theory, however, it is convenient, following [AS], to introduce the mod 2 graded space $\hat{\H}=\H\oplus\H=\H\otimes\C^2$ and to replace Fred$'(\H)$ by Fred$''(\hat{\H})$, the bundle whose fibres are the pairs $(\hat{A},\hat{B})$ of self-adjoint degree 1 operators in $\hat{\H}$ such that $\hat{A}\hat{B}$ and $\hat{B}\hat{A}$ differ from the identity by compact operators.  The space Fred$''(\hat{\H})$ is, of course, homeomorphic to Fred$'(\H)$, but it allows us to use a slightly larger class of twistings.  For if $\H=\H^+\oplus\H^-$ has a mod 2 grading we can give $\hat{\H}=\H\otimes\C^2$ the usual tensor product grading.  As the space Fred$''(\hat{\H})$ of self-adjoint degree 1 operators in $\hat{\H}$ does not change if the grading of $\H$ is reversed, the bundle Fred$''(\hat{P})$ associated to a projective bundle $P$ with involution is well-defined.  It will be technically more convenient, however, to modify the fibre Fred$''(\hat{\H})$ still further, without changing its homotopy type.  Let us recall that for any bounded operator $A$ there is a unique positive self-adjoint operator $\vert A\vert$ such that $\vert A\vert^2=A^*A$.  If now
	\[\hat{A}=\left(\begin{array}{cc} 0 & A\\ A^* & 0\end{array}\right) \ \ \mbox{and} \ \ \hat{B}=\left(\begin{array}{cc} 0 & B^*\\ B & 0\end{array}\right)\]
are self-adjoint degree 1 Fredholm operators which are inverse modulo compact operators then
	\[\tilde{A}=\left(\begin{array}{cc} 0 & \vert B\vert A\\ A^*\vert B\vert & 0\end{array}\right)\]
is another operator of the same type, but with the property that $\tilde{A}^2$ differs from the identity by a compact operator.  It can be connected to $\hat{A}$ in Fred$''(\H)$ by the path $\{\tilde{A}_t\}_{t\in [0,1]}$ where
	\[\tilde{A}_t=\left(\begin{array}{cc} 0 & \vert B\vert^tA\\ A^*\vert B\vert^t & 0\end{array}\right).\]

\begin{defn} If $\hat{\H}$ is a mod 2 graded Hilbert space, let Fred$^{(0)}(\hat{\H})$ denote the space of self-adjoint degree 1 Fredholm operators $\tilde{A}$ in $\hat{H}$ such that $\tilde{A}^2$ differs from the identity by a compact operator, with the topology coming from its embedding $\tilde{A}\mapsto (\tilde{A},\tilde{A}^2-1)$ in ${\mathcal B}\times{\mathcal K}$. \end{defn}

Of course Fred$^{(0)}(\hat{\H})$ is a representing space for $K$-theory, and whenever we have a projective Hilbert bundle $P$ with involution we can define an associated bundle Fred$^{(0)}(P)$.

\begin{defn}
For a projective Hilbert bundle $P$ with involution, we write $K^0_P(X)$ for the space of homotopy classes of sections of Fred$^{(0)}(\hat{P})$, where $\hat{P}=P\otimes {\mathbb P}(\hat{\H})$, where $\hat{\H}$ is a fixed standard mod 2 graded Hilbert space such that both $\hat{\H}^+$ and $\hat{H}^-$ are infinite dimensional.
\end{defn}

Addition in $K^0_P(X)$ is defined by the operation of fibrewise direct sum, so that the sum of two elements naturally lies in $K^0_{P\otimes\P(\C^2)}(X)$, which is canonically isomorphic to $K^0_P(X)$ (see below).  Of course in Fred$'(\H)$ we can define the sum ``internally" simply by composition of operators, but nothing real is gained by that as one needs to pass to $\H\oplus\H$ to see that composition is homotopy-commutative.

\subsubsection*{Remarks}

(i) If $P$ admits a norm-topology structure then Definition 3.3 agrees with the ``naive" definition in terms of sections of Fred$(P)_{\rm norm}$, for the map of bundles 
\[{\rm Fred}(P)_{\rm norm} \rightarrow {\rm Fred}^{(0)}(\hat{P})\]
is a fibre-homotopy equivalence (see [D]).
\bigskip

(ii) The group $K^0_P(X)$ is functorially associated to the pair $(X,P)$, and an isomorphism $\theta:P\rightarrow P'$ of projective bundles induces an isomorphism $\theta_*:K^0_P(X)\rightarrow K^0_{P'}(X)$.  In particular the group Aut$(P)\cong H^2(X;\Z)$ acts naturally on $K^0_P(X)$.  The choice of a definite bundle $P$ representing a class in $H^3(X;\Z)$ is analogous to the choice of a base-point $x_0$ in defining the homotopy group $\pi_i(X,x_0)$, when a path $\gamma$ from $x_0$ to $x_1$ induces
	\[\gamma_*:\pi_i(X,x_0)\stackrel{\cong}{\rightarrow}\pi_i(X,x_1),\]
and $\pi_1(X,x_0)$ acts on $\pi_i(X,x_0)$.  If we give only the class of $P$ in $H^3(X;\Z)$ then the twisted $K$-group is defined only up to the action of $H^2(X;\Z)$. Note, however, that to identify $K^0_{P\otimes\P(\C^2)}$ with $K^0_P$ above we have only to choose an isomorphism between $\hat{\H}\otimes \C^2$ and $\hat \H$, and the space of these isomorphisms is contractible.
\bigskip 

(iii) The standard proof that Fred$(\H)$ is a representing space for $K$-theory (see Appendix 2 or the appendix to[A]) proceeds by showing that a family of Fredholm operators parametrized by a space $X$ can be deformed to a family for which the kernels and cokernels of the operators have locally constant dimension.  These finite dimensional spaces then form vector bundles on $X$, and their difference is the element of $K^0(X)$ corresponding to the family.  In the twisted case, however, such a deformation is never possible if the class of the bundle $P$ in $H^3(X;\Z)$ is not of finite order, for if it were possible then the kernels would define a finite dimensional sub-projective-bundle $P_0$ of $P$, and by Proposition (2.1) (iv) and (vi) the class $[P]=[P_0]$ would have finite order.
\bigskip

(iv) Another peculiarity of twisted $K$-theory when the class $[P]$ is of infinite order is that the index map $K^0_P(X) \rightarrow \Z $ is zero. In other words, any section of Fred$(P)$ takes values in the index zero component of the fibre Fred$(\H)$. This follows easily from the cohomology spectral sequence of the fibration Fred$(P)$ on $X$, a topic which will be examined in our subsequent paper. In particular we shall show that, for the component Fred$_k(P)$ formed by the index $k$ components of the fibres, we have
\[d_3(c_1)=k[P]\]
where $c_1 \in H^2({\rm Fred}(\H );\Z)$ is the universal first Chern class. The spectral sequence gives rise to an exact sequence
\[H^2({\rm Fred}(\H );\Z) \stackrel{d_3}{\rightarrow} H^3(X;\Z) \stackrel{\pi ^*}{\rightarrow} H^3({\rm Fred}_k(P );\Z),\]
where $\pi$ is the projection of the fibre bundle. Thus $\pi ^* \circ d_3 = 0$, and hence $\pi ^* (k[P])=0$. If a section of Fred$_k(P)$ exists then $\pi ^*$ is injective, and hence $k[P]=0$. Since $[P]$ is assumed not to have finite order this implies that $k=0$, as asserted.

\subsubsection*{Algebraic $K$-theory}

We shall now explain how the twisted $K$-theory of a compact space can be defined as the algebraic $K$-theory of a Banach algebra, just as the usual group $K^0(X)$ is the algebraic $K$-theory of the algebra $C(X)$ of continuous complex-valued functions on $X$.  We shall content ourselves with the basic case of twisting by a projective bundle, ignoring bundles with involution.

A bundle $P$ of projective spaces on $X$ gives us a bundle End$(P)$ of algebras, and we might guess that $K^0_P(X)$ is the algebraic $K$-theory of the algebra $\Gamma\,{\rm End}(P)$ of sections of End$(P)$.  This is wrong, however --- even ignoring the problem of topology we encountered in defining Fred$(\H)$ --- unless $P$ is finite dimensional.  If $X$ is an infinite dimensional Hilbert space then $\H\cong\H\oplus\H$, so
	\[{\rm End}(\H)\cong\,{\rm Hom}(\H\oplus\H;\H)\cong\,{\rm End}(\H)\oplus\,{\rm End}(\H)\]
as left-modules over End$(\H)$, and so the algebraic $K$-theory of End$(\H)$ is trivially zero.  Instead of End$(\H)$ we need the Banach algebra ${\mathcal K}=\,{\rm End}_{\rm cpt}(\H)$ of compact operators in $\H$, with the norm topology, which is an algebra without a 1.  The $K$-theory of such a non-unital algebra ${\mathcal K}$ is defined by
	\[K_0({\mathcal K})=\,\mbox{ker}\,:K_0(\hat{\mathcal K})\rightarrow K_0(\C),\]
where $\hat{\mathcal K}=\C\oplus{\mathcal K}$ is the algebra obtained by adjoining a unit to ${\mathcal K}$.  The unital algebra $\hat{\mathcal K}$ has two obvious finitely generated projective modules: $\hat{\mathcal K}$ itself, and also $\H$.  In fact (see [HR])
	\[K_0(\hat{\mathcal K})\cong\Z\oplus\Z\]
with these two generators, and $K_0({\mathcal K})\cong\Z$ with generator $\H$.  (Notice that $\C\otimes_{\hat{\mathcal K}}\H=0$, so $\H$ maps to zero in $K_0 (\C)$.)

With this in mind, we associate to the projective space bundle $P$ the bundle ${\mathcal K}_P$ of non-unital algebras whose fibre at $x$ is End$_{\rm cpt}(P_x)$.  This makes sense because $U(\H)_{c.o.}$ acts continuously on $\K$ (see Appendix 1).

\begin{defn}  The group $K^0_P(X)$ is canonically isomorphic to the algebraic $K$-theory of the Banach algebra $\Gamma (\K_P)$ of sections of $\K_P$. \end{defn}

\noindent {\bf Proof.}  There does not seem to be an obvious map between the two groups, so we shall proceed indirectly, using Bott periodicity ([W],[C],[HR]) for the Banach algebra $\hat{\Gamma}$ formed by adjoining a unit to $\Gamma =\Gamma (\K_P)$.  For $\hat{\Gamma}$, periodicity asserts that $K_0(\hat{\Gamma})\cong\pi_2(BGL(\hat{\Gamma}))\cong\pi_1(GL(\hat{\Gamma}))$, where
	\[GL(\hat{\Gamma})=\bigcup GL_n(\hat{\Gamma})\]
is the infinite general linear group.  We readily deduce
	\[K_0(\Gamma)\cong\lim_{n}\pi_1(GL_n(\Gamma)),\]
where $GL_n(\Gamma)$ denotes the group of invertible $n\times n$ matrices of the form $1+A$, where $A$ has entries in $\Gamma$.  Now $GL_n(\Gamma)$ is the group of sections of the bundle on $X$ associated to $P$ with fibre $GL_n(\K)$.  Furthermore $GL_1(\K)$ is isomorphic to $GL_n(\K)$, and the inclusion
	\[GL_1(\K)\rightarrow GL_n(\K)\]
is a homotopy equivalence.  Finally, $GL_1(\K)$ is known [P] to have the homotopy type of the infinite unitary group $\displaystyle\lim_\rightarrow U_n$, so that its loop-space is $\Z\times BU$.  Putting everything together we find that $K_0(\Gamma)$ is the homotopy classes of sections of the bundle associated to $P$ with fibre $\Z\times BU$, and this is precisely $K^0_P(X)$.
\bigskip

\noindent {\bf Remark.}  The fact that elements of $K^0_P(X)$ cannot be represented by families of Fredholm operators with kernels and cokernels of locally constant dimension corresponds to the existence of two kinds of projective module for $\hat{\K}$ --- ``big" modules like $\hat{\K}$ and ``small" modules like $\H$.  Elements of $K_0(\K)$ can be described using only ``small" modules, but, when we have a twisted family, elements of $K_0(\Gamma\K_P)$ cannot.

\subsection*{More general twistings}

From the point of view of generalized cohomology theories the twistings of $K$-theory which we consider are not the most general possible.  A cohomology theory $h^*$ is represented by a spectrum
	\[h^q(X)\cong [X;h_q],\]
where [\ , \ ] denotes homotopy classes of maps, and $\{h_q\}$ is a sequence of spaces with base-point equipped with homotopy equivalences $h_q\rightarrow\Omega h_{q+1}$.  (Here $\Omega$ denotes the based loop-space.)  Any theory possesses a topological group ${\mathcal G}_h$ of automorphisms which is well-defined up to homotopy.  (In principle an automorphism is a sequence of maps $T_q:h_q\rightarrow h_q$ which commute with the structural maps; but the details of the theory of spectra need great care.)  In any event, the homotopy groups of ${\mathcal G}_h$ are unproblematic: $\pi_i({\mathcal G}_h)$ is the group of transformations of cohomology theories $h^*\rightarrow h^*$ which lower degree by $i$.  Thus if $h^*$ is classical cohomology with integer coefficients ${\mathcal G}_h$ is (up to homotopy) the discrete group $\{\pm 1\}$ of units of $\Z$, for there are no degree-lowering operations.  On the other hand, if $h^*$ is complex $K$-theory then ${\mathcal G}_h$ is much larger.

Whenever we have a principal ${\mathcal G}_h$-bundle $P$ on $X$ we can form the associated bundle of spectra, and can define twisted cohomology groups $h^*_P(X)$.  But for a multiplicative theory $h^*$ --- such as $K$-theory --- it may be natural to restrict to module-like twistings, i.e. those such that $h^*_P(X)$ is a module over $h^*(X)$.  These correspond to a subgroup ${\mathcal G}^{\rm mod}_h$ of ${\mathcal G}_h$ of $\G_h$ with
	\[\pi_0(\G_h^{\rm mod})=h^0(\rm point)^\times\]
	\[\pi_i(\G_h^{\rm mod})=h^{-i}({\rm point}) \ \mbox{for $i>0$}.\]
It is twistings of this kind with which we are concerned here.  We can think of $\G^{\rm mod}_K$ as the ``group" Fred$_{\pm 1}$ of Fredholm operators of index $\pm 1$ under tensor product: it fits into an exact sequence
	\[{\rm Fred}_1\rightarrow\,{\rm Fred}_{\pm 1}\rightarrow (\pm 1).\]
The group Fred$_1$ is a product
	\[{\rm Fred}_1\simeq \P^\infty_{\C}\times S{\rm Fred}_1,\]
where SFred$_1$ is the fibre of the determinant map
\[ {\rm Fred}_1 \cong BU \rightarrow B\T \cong \P^{\infty}_{\C},\]
and the twistings of this paper are those coming from $(\pm 1)\times\P^\infty_{\C}$.  We do not know any equally geometrical approach to the more general ones.

\section{Basic properties of twisted $K$-theory}

In this section we could without any loss use the norm topology on the spaces of Fredholm operators.
\bigskip

One advantage of using the mod 2 graded version Fred$^{(0)}(\hat{P})$ of the bundle of Fredholm operators associated to a projective bundle $P$ is that it gives us at once a multiplication \renewcommand{\theequation}{4.\arabic{equation}}\setcounter{equation}{0}
	\be K^0_P(X)\times K^0_{P'}(X)\rightarrow K^0_{P\otimes P'}(X) \ee
coming from the map
	\[(A,A')\mapsto A\otimes 1+1\otimes A'\]
defined on the spaces of degree 1 self-adjoint Fredholm operators.  (The operator $B = A\otimes 1 + 1\otimes A'$ is Fredholm, as $B^2 = A^2\otimes 1 + 1\otimes A'^2$ because $A\otimes 1$ and $1\otimes A'$ anticommute by the usual conventions of graded algebra, and the positive self-adjoint operator $B^2$ obviously has finite dimensional kernel. If we use the compact-open topology we need to observe that $B^2$ nevertheless varies continuously in the norm topology, so that $\lambda f(\lambda B^2)B$ is a parametrix for $B$ for sufficiently large $\lambda$, where $f:\R \rightarrow \R$ is a smooth function such that $f(t)=t^{-1}$ for $t \geq 1$. We thank J.-L. Tu for pointing out a mistake at this point in an earlier version of this paper.)

In particular, each group $K^0_P(X)$ is a module over the untwisted group $K^0(X)$: this action extends the action of the Picard group Aut$(P)=H^2(X;\Z)$, which is a multiplicative subgroup of $K^0(X)$.  The bilinearity, associativity, and commutativity of the multiplications (4.1) are proved just as for untwisted $K$-theory.

The next task is to define groups $K^i_P(X)$ for all $i\in\Z$, and to check that they form a cohomology theory on the category of spaces equipped with a projective bundle.

The bundle Fred$^{(0)}(\hat{P})$ has a base-point in each fibre, represented by a chosen fibrewise identification $\hat{P}^+_x\cong\hat{P}^-_x$.  We can therefore form the fibrewise iterated loop-space $\Omega^n_X$ Fred$^{(0)}(\hat{P})$, whose fibre at $x$ is $\Omega^n$Fred$^{(0)}(\hat{P}_x)$.  The homotopy-classes of sections of this bundle will be denoted $K^{-n}_P(X)$.  Just as in ordinary $K$-theory these groups are periodic in $n$ with period 2, and we can use this periodicity to define them for all $n\in\Z$.  We have only to be careful to use a proof of periodicity which works fibrewise, i.e. we need a homotopy equivalence
	\[{\rm Fred}^{(0)}(\H)\rightarrow\Omega^2{\rm Fred}^{(0)}(\H)\]
which is equivariant with respect to $U(\H)_{\rm c.o.}$.  The easiest choice is the method of [AS].  For any $n$ we consider the complexified Clifford algebra $C_n$ of the vector space $\R^n$ with its usual inner product.  This is a mod 2 graded algebra, for which we choose an irreducible graded module $S_n$.  Then $S_n\otimes\H$ is also a graded module for $C_n$, and we define Fred$^{(n)}(\H)$ as the subspace of Fred$^{(0)}(S_n\otimes\H)$ consisting operators which commute with the action of $C_n$, in the graded sense.  In [AS] there is defined an explicit homotopy equivalence
	\be {\rm Fred}^{(n)}(S_n\otimes\H)\rightarrow\Omega^n{\rm Fred}^{(0)}(S_n\otimes\H)\cong\Omega^n{\rm Fred}^{(0)}(\H).\ee
On the other hand, when $n$ is even, say $n=2m$, the algebra $C_n$ is simply the full matrix algebra of endomorphisms of the vector space $S_n\cong\C^{2^m}$, and so tensoring with $S_n$ is an isomorphism
	\be {\rm Fred}^{(0)}(\H)\rightarrow\,{\rm Fred}^{(n)}(S_n\otimes\H).\ee
The maps (4.2) and (4.3) are completely natural in $\H$, and make sense fibrewise in Fred$^{(0)}(P)$.

To be a cohomology theory on spaces with a projective bundle means that $K^*_P$ must be homotopy-invariant and must possess the Mayer-Vietoris property that if $X$ is the union of two subsets $X_1$ and $X_2$ whose interiors cover $X$, and $P$ is a projective bundle on $X$, there is an exact sequence
	\[\ldots \stackrel{d}{\longrightarrow}K^i_P(X)\rightarrow K^i_{P_1}(X_1)\oplus K^i_{P_2}(X_2)\rightarrow K^i_{P_{12}}(X_{12})\stackrel{d}{\longrightarrow}K^{i+1}_P(X)\rightarrow\ldots\]
where $X_{12}=X_1\cap X_2$, and $P_1,P_2,P_{12}$ are the restrictions of $P$ to $X_1,X_2,X_{12}$.  The proof of this is completely standard, and we shall say no more about it than that the definition of the boundary map $d$, when $i=-1$, is as follows.  One chooses $\varphi :X\rightarrow [0,1]$ such that $\varphi\vert X_1=0$ and $\varphi\vert X_2=1$.  Then if $s$ is a section of $\Omega_X$Fred$^{(0)}(P)$ defined over $X_{12}$ we define the section $ds$ of Fred$^{(0)}(P)$ to be the base-point outside $X_{12}$, and at $x\in X_{12}$ to be the evaluation of the loop $s(x)$ at time $\varphi(x)$.

\subsubsection*{The spectral sequence}

Once we have a cohomology theory we automatically have a spectral sequence defined for any space $X$ with a projective bundle $P$, relating $K^*_P(X)$ to classical cohomology.  More precisely,

\begin{prop} There is a spectral sequence whose abutment is $K^*_P(X)$ with
	\[E^{pq}_2=H^p(X;K^q({\rm point})).\]
The coefficients here are twisted by the class $\xi_P$ of $P$ in $H^1(X;\Z/2)$. \end{prop}

  The spectral sequence is constructed exactly as in the untwisted case, e.g. by the method of [S1].  We shall discuss this further in the sequel to this paper, where we shall determine the first non-zero differential $d_3$, and shall use the spectral sequence to describe $K^*_P(X)\otimes \mathbb Q$.

\section{Examples}

An important source of projective spaces which do not have canonically defined underlying vector spaces is the fermionic Fock space construction, due originally to Dirac.  If $\H$ is a Hilbert space with an orthonormal basis $\{e_n\}_{n\in\Z}$ one can consider the Hilbert space $\F(\H)$ spanned by an orthonormal basis consisting of the formal symbols
	\[e_{n_1}\wedge e_{n_2}\wedge e_{n_3}\wedge \ldots\]
where $n_1>n_2>n_3>\ldots$ and $n_{k+1}=n_k-1$ for all large $k$.  We can think of $\F(\H)$ as a ``renormalized" version of the exterior algebra of $\H$.  The important thing for our purposes is that the projective space $\P\F(\H)$ of $\F(\H)$ depends only weakly on the choice of the orthonormal basis $\{e_n\}$.  Because
\[ \F (\H )\cong \Lambda (\H ^+ ) \otimes \Lambda (\bar{\H}^- )\]
  it clearly depends only on the decomposition $\H=\H^+\oplus\H^-$, where $\H^+$ is spanned by $\{e_n\}_{n\geq 0}$;  but, less obviously, it depends only on the {\it polarization} of $\H$, i.e. on the {\it class} of the decomposition in a sense explained in [PS] Chap. 7.  The case of interest here is when $\H=\H_E$ is the space of sections of a smooth complex vector bundle $E$ on an oriented circle $S$.  If we choose a parametrization $\theta :S\rightarrow\R/2\pi\Z$ and a trivialization $E\cong S\times\C^m$ then the class of  the splitting for which $\H^+$ is spanned by $v_ke^{in\theta}$ for $n\geq 0$, where $\{v_k\}$ is the basis of $\C^m$, is independent of both the parametrization $\theta$ and the trivialization, so that the projective space
	\[\P_E=\P\F(\H_E)\]
depends only on $E$.

We can apply this as follows.  For each element $u$ of the unitary group $U_m$ let $E^u$ be the vector bundle on $S^1=\R/2\pi\Z$ with holonomy $u$. (In other words, $E^u$ is obtained from $\R \times \C^n$ by identifying $(x + 2\pi ,\xi )$ with $(x,u\xi)$.)  Then the spaces $\P_{E^u}$ form a projective bundle on the group $U_m$.  We shall denote this bundle again by $\P_E$: we hope the notation will not prove confusing.  The bundle $\P_E$ on $U_m$ is equivariant with respect to the action of $U_m$ on itself by conjugation: an element $g\in U_m$ defines an isomorphism $E^u\rightarrow E^{gug^{-1}}$, and hence an isomorphism $\P_{E^u}\rightarrow \P_{E^{gug^{-1}}}$.  We shall return to this aspect of the bundle in \S 6.  We can also regard $\P_{E}$ as a projective bundle with involution, for multiplication by $\{\pm 1\}$ on $\H$ induces a projective action of the group $\{\pm 1\}$ on $\F(\H)$.

\begin{prop} The class of the projective bundle $\P_E$ on $U_m$ is a generator of $H^3(U_m;\Z)\cong\Z$, and as a bundle with involution its class is the non-zero element of $H^1(U_m;\Z/2)\cong\Z/2$.\end{prop}

Before justifying this assertion we shall mention a similar example, which is actually the one used by Freed, Hopkins, and Teleman.  For a finite dimensional complex vector space $W$ with an inner product the projective space of the exterior algebra $\Lambda (W)$ is independent of the complex structure on $W$, as it is canonically isomorphic to the projective space of the spin module $\Delta (V)$ of the real vector space $V$ underlying $W$.  Another way of saying this is that if we start with an even-dimensional real vector space $V$ then there is a canonical factorization of complex projective spaces
\renewcommand{\theequation}{5.\arabic{equation}}\setcounter{equation}{1}
	\be \P(\Lambda (V_{\C}))\cong\P(\Delta (V))\otimes\P(\Delta (V)), \ee
where $V_{\C}$ is the complexification of $V$.  There is an infinite dimensional analogue of this phenomenon, explained in Chapter 12 of [PS].  If $\H$ is a real Hilbert space a {\it complex polarization} of $\H$ will mean a preferred class of complex structures --- equivalently, a class of decompositions $\H_{\C}=\H^+\oplus\H^-$ with $\H^+$ and $\H^-$ complex conjugate.  If $\H$ has a complex polarization then we can define a projective spin module $\P (\Delta (\H))$, and
	\begin{eqnarray} \P\F(\H_{\C})&\cong&\P(\Lambda(\H^+)\otimes\Lambda(\bar{\H}^-))\nonumber\\ &\cong&\P(\Delta (\H))\otimes\P(\Delta (\H)).\end{eqnarray}

	Before applying this to bundles on the circle we need a little more discussion.  The first point is that the isomorphisms (5.2) and (5.3) are functorial in the category of projective spaces {\it with involution}.  This is important because an orientation-reversing automorphism of $V$ interchanges the components of $\Delta (V)$.  Next, if we have an {\it odd}-dimensional real vector space $V$ we define $\Delta (V)=\Delta (V\oplus\R)$, but we must think of it as having an additional action of the Clifford algebra $C_1$ on one generator (commuting in the graded sense with the action of the Clifford algebra $C(V)$ which $\Delta (V)$ possesses in all cases).  For odd dimensional $V$ the isomorphism (5.2) is replaced by
	\[\P(\Lambda (V_{\C}))\otimes\P(S_2)\cong\P(\Delta(V))\otimes\P(\Delta (V)),\]
as projective spaces with involution, where, on the left, the space $S_2\cong\C^2$ is the irreducible module for the Clifford algebra $C_2\cong C_1\otimes C_1$.  There is exactly the same distinction between ``odd" and ``even" dimensionality for polarized real Hilbert spaces $\H$, according as $\H$ or $\H\oplus\R$ has a preferred class of complex structures.

Now let us consider the real Hilbert space $\H_E$ of sections of a smooth real vector bundle $E$ on the circle $S^1$.  The Fourier decomposition gives either $\H_E$ or $\H_E\oplus\R$ a class of complex structures: in fact $\H_E$ is ``even-dimensional" if $E$ is even-dimensional and orientable, or if $E$ is odd-dimensional and non-orientable, and $\H_E$ is ``odd-dimensional" otherwise.  We shall write $\P_E^{\rm spin}$ for the projective Hilbert space $\P\Delta(\H_E)$.  As before, we can consider the family of $m$-dimensional real bundles $E^u$ on $S^1$ parametrized by elements $u$ of the orthogonal group $O_m$.  The corresponding projective spaces $\P_{E^u}^{\rm spin}$ form a bundle $\P_E^{\rm spin}$ on $O_m$.

\setcounter{prop}{3}
\begin{prop} The class of the bundle $\P_E^{\rm spin}$ --- with its involution --- on $O_m$ is $(\varepsilon ,\eta)\in H^1(O_m:\Z/2)\oplus H^3(O_m;\Z)$, where $\varepsilon$ restricts to the non-trivial element, and $\eta$ to a generator, on each connected component of $O_m$. \end{prop}

To prove Propositions (5.1) and (5.4), let us take a slightly different point of view on the preceding constructions.  If $G$ is a compact connected Lie group, let ${\mathcal L} G$ denote the group of smooth loops $S^1=\R/2\pi\Z\rightarrow G$, and let ${\mathcal P}G$ be the space of smooth maps $f:\R\rightarrow G$ such that $\theta\mapsto f(\theta +2\pi)f(\theta)^{-1}$ is constant.  Then ${\mathcal L}G$ acts freely on ${\mathcal P}G$ by right multiplication, and the map ${\mathcal P}G\rightarrow G$ given by $f\mapsto f(2\pi)f(0)^{-1}$ makes ${\mathcal P}G$ a principal ${\mathcal L}G$-bundle over $G$.  Thus for any projective representation $\P$ of ${\mathcal L}G$ we have an associated projective bundle ${\mathcal P}G\times_{{\mathcal L}G}\P$ on $G$ --- in fact a $G$-equivariant bundle, when $G$ acts on itself by conjugation, in view of the action of $G$ on ${\mathcal P}G$ by left multiplication.  The invariant of ${\mathcal P}G\times_{{\mathcal L}G}\P$ in $H^3(G;\Z)$ is clearly represented by the composite
	\[G\rightarrow B{\mathcal L}G\rightarrow BPU(\H)\simeq K(\Z,3),\]
where the first map is the classifying map for ${\mathcal P}G$ and the second is induced by the representation ${\mathcal L}G\rightarrow PU(\H)$.  This implies that the transgression $H^3(G;\Z)\rightarrow H^2({\mathcal L}G;\Z)$ takes the invariant to the class of the circle bundle on ${\mathcal L}G$ which is the central extension defined by $\P$.  The bundle $\P_E$ on $U_m$ which we described above is obtained from ${\mathcal P}U_m$ by what is called the {\it basic representation} of ${\mathcal L}U_m$.  (To see this, think of an element of ${\mathcal P}U_m$ over $u\in U_m$ as defining a trivialization of the bundle $E^u$.)  Because the maps
	\[H^3(U_m;\Z)\rightarrow H^3(SU_m;\Z)\rightarrow H^2({\mathcal L}SU_m;\Z)\cong\Z\]
are isomorphisms, we need only ask which central extension of ${\mathcal L}SU_m$ acts on the basic representation, and we know from [PS] that we get a generator of $H^2({\mathcal L}SU_m;\Z)$.  The other part of (5.1), concerning the class in $H^1(U_m;\Z/2)$, is much easier, as all we need to know is that an element of ${\mathcal L}U_m$ of winding number 1 acts on the Fock space ${\mathcal F}(L^2(S^1;\C^n))$ by an operator which raises degree by 1.

Proposition (5.4) follows easily from (5.1).  First, one may as well assume $m$ is even.  Then the bundle $\P_E^{\rm spin}$ on $O_{2k}$ restricts to $\P_E$ on $U_k$, while the maps
	\[H^1(SO_{2k};\Z/2)\rightarrow H^1(U_k;\Z/2)\]
and
	\[H^3(SO_{2k};\Z)\rightarrow H^3(U_k;\Z)\]
are isomorphisms; this deals with the invariants on the identity component of $O_{2k}$.  The other component can be treated by embedding $U_{k-1}$ in it by adding a fixed non-orientable bundle and using the multiplicativity of the Fock space construction.
\bigskip

Let us now describe some families of Fredholm operators in the projective bundles we have just constructed.  In the representation theory of a loop group ${\mathcal L}G$ one usually studies projective representations $\H$ which are of {\it positive energy} and {\it finite type}.  This means that the circle ${\mathbb T}$ of rotations of the loops acts unitarily on $\H$, compatibly with its action on ${\mathcal L}G$, and decomposes $\H$ into finite dimensional eigenspaces
	\[\H=\bigoplus_{n\geq 0}\H_n,\]
where ${\mathbb T}$ acts on $\H_n$ by the character $e^{i\theta}\mapsto e^{in\theta}$.   (One calls $\H_n$ the part of ``energy" $n$.)  The infinitesimal generator $L_0$ of the circle action is an unbounded positive self-adjoint operator in $\H$.  When we consider the family ${\mathcal P}\times_{{\mathcal L}G}\P(\H)$ on $G$ the group $\R$ acts on ${\mathcal P}$ by translation, compatibly with the action of ${\mathbb T}=\R/2\pi\Z$ on ${\mathcal L}G$ and $\P(\H)$.  So $\R$ acts fibrewise on the bundle.  If we identify the fibre ${\mathbb P}_g$ at $g\in G$ with ${\mathbb P}(\H)$ by choosing $f\in{\mathcal P}$ such that $f(\theta +2\pi)f(\theta)^{-1}=g$ then the infinitesimal generator $L_0^{(g)}$ of the $\R$-action on $\P_g$ is clearly given by
	\[L_0^{(g)}=L_0+f^{-1}f',\]
where $f^{-1}f'$, which is periodic, is regarded as an element of the Lie algebra of ${\mathcal L}G$.  In fact we can choose $f$ to be a 1-parameter subgroup of $G$ generated by an element $\xi\in{\mathfrak g}=\,{\rm Lie}(G)$ such that exp$(2\pi\xi)=g$, and then
	\[L_0^{(g)}=L_0+\xi.\]
As $\xi$ commutes with $L_0$ it acts separately in each energy level $\H_n$.  In fact we know from [PS](9.3.7) that if $V_\lambda$ is an irreducible representation of $G$ with highest weight $\lambda$ contained in $\H_n$ then $\|\lambda\|^2 \leq an+b$, where $a$ and $b$ are constants depending on the representation $\H$. On the other hand the eigenvalues of $\xi$ in $V_{\lambda}$ are bounded by $\|\lambda\|\|\xi \|$, so the eigenvalues of $\xi$ in $\H _n$ grow only like $n^{1/2}$ as $n \rightarrow \infty$. This shows that for any $g \in G$ the operator $L_0^{(g)}$ decomposes the Hilbert space $\H_g$ underlying the projective space $\P_g$ into the orthogonal direct sum of a sequence of finite-dimensional eigenspaces $\H_{g,\lambda}$ corresponding to a sequence of eigenvalues $\{\lambda \}$, depending on $g$ and tending to $\infty$. In particular, the zero-eigenspace of $L_0^{(g)}$ is always finite-dimensional.

\bigskip

The family $\{L_0^{(g)}\}$, being positive, is not directly interesting in $K$-theory. It is analogous to the family of Laplace operators on the fibres of a bundle of compact manifolds, and we need something analogous to the family of Dirac operators.   For a positive energy representation $\H$ of a loop group $\L G$  Freed, Hopkins, and Teleman consider the projective bundle $\P(\H)_G={\mathcal P}G\times_{\L G}\P(\H)$ on $G$ which we have already described.  Its fibre $\P_g =\P (\H_g)$ at $g\in G$ is a representation of the twisted loop group $\L_gG$ whose Lie algebra $\L_g \mathfrak g$ is the space of sections of the real vector bundle $E_g$ on $S^1$ with fibre ${\mathfrak g}$ and holonomy $g$.  They tensor $\P(\H)_G$ with the spinor bundle $\P^{\rm spin}_E$.  There is then a Dirac-type operator $D_{\H}= \{D_g \}$ acting fibrewise in $\P(\H)_G\otimes\P_E^{\rm spin}$, defined for $\xi \otimes \psi \in \H_g \otimes \Delta(\L_g \mathfrak g )$ by
\[D_g (\xi \otimes \psi )= \sum e_i \xi \otimes e_i^* \psi ,\]
where $\{e_i \}$ is a basis of $L_g \mathfrak g ^*$, and $\{e^*_i \}$ is the dual basis of $\L_g \mathfrak g ^*$, regarded as elements of the Clifford algebra $C(\L_g \mathfrak g ^* )$. (If $\xi$ and $\psi$ are in $L_0^{(g)}$-eigenspaces, and we choose the basis $\{e_i \}$ to consist of $L_0^{(g)}$-eigenvectors in $\L_{(g)} \mathfrak g _{\C}$, then the sum on the right is finite.) The operator $D_g$ is, of course, an unbounded operator, but of a very tractable kind. It is defined on the dense subspace which is the algebraic direct sum of the finite-dimensional eigenspaces of $L_0^{(g)}$, and its square is a scalar multiple of $L_0^{(g)}$. It therefore decomposes as the sum of finite-dimensional operators acting in the $L_0^{(g)}$-eigenspaces. We can obtain a family $\{ A_g \}$ of bounded Fredholm operators from the family $ \{D_g \}$ by defining 
\[A_g = (D_g^2 +1)^{-1/2}D_g.\]
The family $\{ A_g \}$ defines an element of the twisted $K$-theory of $G$ --- in fact of the $G$-equivariant twisted $K$-theory --- for each projective representation $\H$ of the loop group $\L G$.  This is the map which Freed-Hopkins-Teleman prove to be an isomorphism.  (If $G$ is odd-dimensional, so is, as we have seen, the polarized Lie algebra $\L{\mathfrak g}$, and then the additional $C_1$-action on $\P_E^{\rm spin}$ gives us an odd-dimensional $K$-theory class.)

\setcounter{section}{5}
\section{The equivariant case}

When a compact group $G$ acts on a space $X$ we can define equivariant $K$-theory $K^*_G(X)$.  If $X$ is compact then $K^0_G(X)$ is the Grothendieck group of $G$-vector-bundles on $X$.  If $X$ is not compact, however, then one normally defines $K^0_G(X)$ as the equivariant homotopy classes of $G$-maps from $X$ to a suitable representing $G$-space $\k^0_G$.  Just as in the non-equivariant case, the space $\k^0_G$ can be chosen in quite a variety of ways.  If $\H_G$ is what we shall call a {\it stable} $G$-Hilbert-space, i.e. a Hilbert space representation of $G$ in which each irreducible representation of $G$ occurs with infinite multiplicity (or, equivalently, one such that $\H_G \cong \H_G \otimes L^2 (G)$), then any $G$-vector-bundle on a compact base-space $X$ can be embedded as a $G$-subbundle of $X\times\H_G$, and so can be pulled back from the Grassmannian Gr$(\H_G)$ of all finite dimensional vector subspaces of $\H_G$.  Stabilizing in a familiar way gives us a natural candidate for $\k^0_G$.  (A convenient choice of the stabilization is the {\it restricted} Grassmannian Gr$_{\rm res}(\H_G)$ mentioned in \S 3.)

The space Fred$(\H_G)$ of Fredholm operators in $\H_G$, with the norm topology, might seem another natural choice for $\k^0_G$, but unfortunately the action of $G$ on Fred$(\H_G)$ is very far from continuous.  This can be dealt with in two ways. One is to replace Fred$(\H_G)$ by the $G$-continuous subspace
	\[{\rm Fred}_{G-{\rm cts}}(\H_G)=\{A\in\,{\rm Fred}(\H_G):g\mapsto gAg^{-1}\ \mbox{is continuous}\},\]
which is closed in Fred$(\H_G)$, and is a representing space for $K^0_G$, as is proved in Appendix 3. The other way is to pass to the more sophisticated space Fred$^{(0)}(\H_G)$ introduced in \S3.

To twist equivariant $K$-theory we need a bundle $P$ of projective spaces on which $G$ acts, mapping $P_x$ to $P_{gx}$ by a projective isomorphism.  We shall call $P$ {\it stable} if $P\cong P\otimes L^2(G)$.   As before, we must decide whether or not to require that the structural group of $P$ is $U({\mathcal H})$ with the norm topology.  Either way, we must be more careful than in the non-equivariant case.  If $P$ has structural group $U({\mathcal H})_{\rm norm}$ when the $G$-action is ignored it is impossible for $G$ to act continuously on the associated principal bundle of $P$ (unless $G$ acts almost freely on $X$).  Instead, we must require that

(i) each point $x\in X$ with isotropy group $G_x$ has a $G_x$-invariant neighbourhood $U_x$ such that there is an isomorphism of bundles with $G_x$-action
	\[P\vert U_x\cong U_x\times {\mathbb P}({\mathcal H}_x)\]
for some projective space ${\mathbb P}({\mathcal H}_x)$ with $G_x$-action, and

(ii) the transitions between these trivializations are given by maps
	\[U_x\cap U_y\rightarrow \,{\rm Isom}({\mathcal H}_x;{\mathcal H}_y)\]
which are continuous in the norm topology.
\bigskip

When $P$ satisfies these conditions we can associate to it the bundle Fred$(P)$, defined without using the $G$-action of $P$, and with the norm topology in each fibre.  Although the natural action of $G$ on Fred$(P)$ is not continuous, it makes sense to define $K^0_{G,P}(X)$ as the group of homotopy classes of $G$-equivariant continuous sections of Fred$(P)$.

As in the non-equivariant case, however, we prefer to avoid the norm topology.  For any locally trivial projective bundle $P$ with $G$-action the group $G$ acts continuously on the associated bundle Fred$^{(0)}(P)$.  Even using Fred$^{(0)}(P)$, however, it seems essential to require the bundle $P$ to satisfy condition (i): otherwise we do not, for example, see how to show that Fred$^{(0)}(P)$ is equivariantly trivial when $P={\mathbb P}(E)$ comes from a stable equivariant bundle $E$ of Hilbert spaces on $X$ (cf. the action of $G=(\pm 1)$ on $E=[0,1]\times L^2([0,1])$ given by
	\[(-1).(x,\phi)=(x,\varepsilon_x\phi),\]
where
	\begin{eqnarray*} \varepsilon_x(y)&=&1 \  \ \mbox{when} \ \ y\leq x\\ &=&-1 \ \ \ \mbox{when} \ \ y>x.)\end{eqnarray*}
If condition (i) holds then we can trivialize $E$ over a compact base $X$ by constructing a $G$-equivariant section of the bundle on $X$ with fibre Isom$({\mathcal H}_G;E_x)$ at $x$.  This can be done by induction on the number of sets in a covering of $X$ by $G$-invariant open sets of the form $G.S_i$, where $S_i$ is a $G_{x_i}$-invariant ``slice" (see [Bor] Chap.7, and [S2] page 144) at a point $x_i\in X$, and $E\vert S_i$ is $G_{x_i}$-equivariantly trivial.

\begin{defn} 
For stable projective bundles $P$ which satisfy condition (i) above we define $K^0_{G,P}(X)$ as the group of homotopy classes of equivariant sections of Fred$^{(0)}(P)$.
\end{defn}

The passage from twisted $K$-theory to the equivariant twisted theory is now quite unproblematical, at any rate for those accustomed to ordinary equivariant $K$-theory [S2].  There seems no point in spelling it out.  The most interesting thing to discuss is the classification of stable $G$-projective-bundles $P$, i.e. the analogue of Proposition (2.1) and Proposition (2.2).  A $G$-projective bundle has an invariant $\eta_P$ in the equivariant cohomology group $H^3_G(X;\Z)$.  This group can be defined by means of the ``Borel construction", i.e. the functor which takes a $G$-space $X$ to $X_G=(X\times EG)/G$, where $EG$ is a fixed contractible space on which $G$ acts freely.  

\begin{defn}
	\[H^*_G(X;\Z)=H^*(X_G;\Z).\]
\end{defn}

In particular, $H^*_G({\rm point};\Z)=H^*(BG;\Z)$, where $BG$ is the classifying space $EG/G$.
\bigskip

Let us write Pic$_G(X)$ for the group of isomorphism classes of complex $G$-line-bundles on $X$ (or, equivalently, of principal $\T$-bundles with $G$-action), and Proj$_G(X)$ for the group of stable $G$-projective-bundles satisfying condition (i). Applying the Borel construction to line bundles and projective bundles gives us homomorphisms
	\[{\rm Pic}_G(X) \rightarrow \,{\rm Pic}(X_G)\cong H^2_G(X;\Z)\]
	\[{\rm Proj}_G(X) \rightarrow \,{\rm Proj}(X_G)\cong H^3_G(X;\Z),\]
which we shall show are bijective.
\bigskip

\no {\bf Remark}
A mod 2 graded projective bundle, in the sense of \S 2, is a projective bundle with $\Z/2$-action on a base $X$ with trivial $\Z/2$-action. If $G=\Z/2$ acts trivially on $X$ then
\[H_G^*(X)=H^* (X\times {\mathbb R} P^{\infty}) \cong H^* (X;H^* ({\mathbb R}P^{\infty})),\]
so that
\[H^3_G(X;\Z) \cong H^1 (X;\Z/2) \oplus H^3 (X;\Z ).\]
This agrees set-theoretically with 2.3, but the tensor product of $G$-spaces is not the same as the graded tensor product.

\begin{prop}
\begin{enumerate}\item[(i)] $H^2_G({\rm point};\Z)\cong\,{\rm Hom}(G;\T)$
\item[(ii)] $H^3_G({\rm point};\Z)\cong\,{\rm Ext}(G;\T)$, the group of central extensions
\[1\rightarrow \T \rightarrow \tilde{G} \rightarrow G \rightarrow 1.\]
\item[(iii)] ${\rm Pic}_G(X)\cong H^2_G(X;\Z)$
\item[(iv)] ${\rm Proj}_G(X)\cong H^3_G(X;\Z)$, and this remains true if we replace the left-hand side by the group of stable $G$-projective bundles with norm-topology structural groups.
\end{enumerate}
\end{prop}

Of course the assertions (i) and (ii) here follow from (iii) and (iv), but they are easier to prove, and seem worth making explicit. Because ${\mathbb P}_G={\mathbb P}(\H _G)$ is a classifying space for $G$-line-bundles when $\H _G$ is an ample $G$-Hilbert-space, (iii) is simply the fact that ${\mathbb P }_G$ represents the functor $H^2_G(\ \ ;\Z)$, which can be proved quite easily in a variety of ways. The method we follow is chosen for its wider applications.
\bigskip

Before giving the proof of 6.3, let us review the bundles of Fock spaces on a group $G$ which were described in \S 5. These bundles are $G$-equivariant when $G$ acts on the base-space by conjugation. They satisfy the equivariant local triviality conditiion (i) because the principal fibration $\p G\rightarrow G$ has the corresponding property. They are not the most general possible equivariant bundles, as the action of the isotropy group on each fibre extends (non-canonically) to an action of $G$. They do not, however, have a natural norm-continuous structure, for the natural identifications of the fibre $\P_g$ at $g$ with $\P_1$ differ among themselves by the action of elements of $\L G$ on $\P_1$, and so the natural transition maps between local trivializations will factorize through $\L G$, which sits as a discrete subspace in $U(\H)_{\rm norm}$.

These equivariant projective bundles are determined by their classes in $H^3_G(G_{\rm conj};\Z)$. The Borel construction $EG \times _G G_{\rm conj}$ is simply the free loop space $\L BG$, which for connected $G$ is the same as $B\L G$. In the connected case this is most clearly seen by writing
\[EG \times _G G_{\rm conj}=EG \times _G (\p/\L G)=(EG \times _G \p)/\L G \simeq B\L G,\]
as $G\backslash \p$ can be identified with the affine space of connections in the trivial $G$-bundle on the circle, so that $EG \times _G \p$ is a contractible space on which $\L G$ acts freely. From this point of view it is clear that the class of the bundle on $G_{\rm conj}$ coming from a projective representation $\H$ of $\L G$ is simply the topological class of the bundle
\[B\T \rightarrow B\tilde{\L} G \rightarrow B\L G\]
with fibre $B\T \simeq \P^\infty _\C$, where $\tilde{\L} G$ is the central extension of $\L G$ which acts on $\H$.

If $G$ is connected and semisimple, the Serre spectral sequence for $H^*_G$ gives us an exact sequence
\[0\rightarrow {\rm Ext}(G;\T )\rightarrow H^3_G (G;\Z )\rightarrow  H^3 (G;\Z ),\]
where the inclusion of Ext$(G;\T )$ is split by restriction to $1 \in G$. Thus the class of an equivariant projective bundle --- or of a representation of $\L G$ --- is determined by its non-equivariant class together with its class as a projective representation of $G$, and the examples of \S 5 show us that when $G= {\rm SO}_m$ any class in $H^3 (G;\Z )$ can arise.

When $G=U_m$, on the other hand, the spectral sequence gives us an exact sequence
\[0\rightarrow H^2 (BU_m ;H^1 (U_m ;\Z )) \rightarrow H^3_{U_m}(U_m ; \Z ) \rightarrow H^3 (U_m ; \Z ) \rightarrow 0. \]
When $m=1$ this tells us that $H^3_{U_1}(U_1 ; \Z ) \cong \Z$, the invariant being the flow of the grading of a $\Z$-graded projective bundle around the base circle. When $m>1$, we have
\[H^3_{U_m}(U_m ; \Z ) \cong \Z \oplus \Z \]
by the map
\[H^3_{U_m}(U_m ; \Z ) \rightarrow H^3_{U_1}(U_1 ; \Z )\oplus H^3 (U_m;\Z ).\]

\bigskip  

To prove (6.3) it is helpful to introduce groups $H^*_G(X;A)$ defined for any topological abelian group $A$.  These are the hypercohomology groups of a simplicial space $\X$ whose ``realization" is the space $X_G$ (see [S1]).  Whenever a group $G$ acts on a space $X$ we have a topological category whose space of objects is $X$ and whose space of morphisms from $x_0$ to $x_1$ is $\{g\in G:gx_0=x_1\}$.  (Thus the complete space of morphisms is $G\times X$.) A topological category can be regarded as a simplicial space $\X$ whose space $\X_p$ of $p$-simplexes is the space of composable $p$-tuples of morphisms in the category: in our case $\X_p=G^p\times X$.

For any simplicial space $\X$ and any topological abelian group $A$ we can define the hypercohomology $\h^*(\X;sh(A))$ with coefficients in the sheaf of continuous $A$-valued functions.  It is the cohomology of a double complex $C^{..}$, where, for each $p\geq 0$, the cochain complex $C^{p.}$ calculates $H^*(\X_p;sh(A))$.

\begin{defn}
	\[H^*_G(X;A)=\h^*(\X;sh(A)).\]
\end{defn}

If $A$ is discrete, the hypercohomology is just a way of calculating the cohomology of the realization $X_G$ of $\X$, so the new definition of $H^*_G(X;A)$ agrees with the old one.  In any case, the groups $H^*_G(X;A)$ are the abutment of a spectral sequence with
	\[E_1^{pq}=H^q(G^p\times X;sh(A)).\]

\begin{lem}
If $G$ is a compact group, then
	\[H^{p+1}_G(X;\Z)\cong H^p_G(X;\T)\]
for any $p>0$.
\end{lem}

\no {\bf Proof} \ \ Because of the exact sequence
	\[0\rightarrow sh(\Z)\rightarrow sh(\R)\rightarrow sh(\T)\rightarrow 0\]
it is enough to show that $H^p_G(\X;\R)=0$ for $p>0$.  As $E^{pq}_1=0$ for $q>0$ in the specctral sequence when $A=\R$, we see that $H^*_G(X;\R)$ is simply the cohomology of the cochain complex of continuous real-valued functions on the simplicial space $\X$, which is easily recognized as the complex of continuous Eilenberg-Maclane cochains of the group $G$ with values in the topological vector space Map$(X;\R)$ of continuous real-valued functions on $X$.  This complex is well-known to be acyclic in degrees $>0$ when $G$ is compact.  (It is the $G$-invariant part of the contractible complex of so-called ``homogeneous cochains", and taking the $G$-invariants is an exact functor, simply because cochains can be averaged over $G$.)
\medskip

\no {\bf Proof of (6.3)} \ \ (i) When $X$ is a point we have $E^{0q}_1=0$ in the spectral sequence for $H^*_G$, and we have already pointed out that $E^{p0}_2=H^p_{\rm c.c.}(G;A)$ is the cohomology of $G$ defined by continuous Eilenberg-Maclane cochains.  So
	\[H^1_G({\rm point};A)\cong E^{10}_2\cong H^1_{\rm c.c.}(G;A)\cong \,{\rm Hom}(G,A)\]
for any topological abelian group $A$.
\bigskip

(ii) In this case the spectral sequence gives us an exact sequence
	\[0\rightarrow E^{20}_2\rightarrow H^2_G({\rm point};\T)\rightarrow E^{11}_2\rightarrow E^{30}_2,\]
i.e.
	\[0\rightarrow H^2_{\rm c.c.}(G;\T)\rightarrow H^2_G({\rm point};\T)\rightarrow\,{\rm Pic}(G)_{\rm prim}\rightarrow H^3_{\rm c.c.}(G;\T),\]
for $E^{11}_1=H'(G;sh(\T))=\,{\rm Pic}(G)$, and $E^{11}_2$ is the subgroup of primitive elements, i.e. of circle bundles $\tilde{G}$ on $G$ such that
	\[m^*\tilde{G}\cong pr^*_1\tilde{G}\otimes pr^*_2\tilde{G},\]
where $pr_1,pr_2,m:G\times G\rightarrow G$ are the obvious maps.  Equivalently, Pic$(G)_{\rm prim}$ consists of circle bundles $\tilde{G}$ on $G$ equipped with bundle maps $\tilde{m}:\tilde{G}\times\tilde{G}\rightarrow\tilde{G}$ covering the multiplication in $G$.  It is easy to see that the composite
\renewcommand{\theequation}{6.\arabic{equation}}\setcounter{equation}{5}
	\be {\rm Ext}(G;\T)\rightarrow H^2_G({\rm point};\T)\rightarrow\,{\rm Pic}(G) \ee
takes an extension to its class as a circle bundle.  On the other hand $H^2_{\rm c.c.}(G;\T)$ is plainly the group of extensions $\T\rightarrow\tilde{G}\rightarrow G$ which as circle bundles admit a continuous section, so its image in Ext$(G;\T)$ is precisely the kernel of (6.6).  It remains only to show that the image of Ext$(G;\T)$ in Pic$(G)_{\rm prim}$ is the kernel of
	\[{\rm Pic}(G)_{\rm prim}\rightarrow H^3_{\rm c.c.}(G;\T).\]
This map, however, associates to a bundle $\tilde{G}$ with a bundle map $\tilde{m}$ as above precisely the obstruction to changing $\tilde{m}$ by a  bundle map $G\times G\rightarrow\T$ to make it an associative product on $\tilde{G}$.
\bigskip

(iii) The spectral sequence gives
	\[0\rightarrow E^{10}_2\rightarrow H^1_G(X;\T)\rightarrow E_2^{01}\rightarrow E^{20}_2.\]
Now $E^{01}_2=P_{ic}(X)$, and $E^{01}_2$ is the subgroup of circle bundles $S\rightarrow X$ which admit a bundle map $\tilde{m}:G\times S\rightarrow S$ covering the $G$-action on $X$.  As before, $\tilde{m}$ can be made into a $G$-action on $S$ if and only if an obstruction in $H^2_{cc}(G;{\rm Map}(X;\T))$ vanishes.  Finally, the kernel of Pic$_G(X)\rightarrow\,{\rm Pic}(X)$ is the group of $G$-actions on $X\times\T$, and this is just $E^{10}_2=H^1_{cc}(G;{\rm Map}(X;\T))$.
\bigskip

(iv) This is the essential statement for us, and is distinctly harder to prove than the other three.  If we knew a priori that the functor $X\mapsto\,{\rm Proj}_G(X)$ was representable by a $G$-space the argument would be much simpler; but we do not see a simple proof of representability.  Instead we shall prove by the preceding methods that the map Proj$_G(X)\rightarrow H^3_G(X;\Z)$ is injective, and then we shall construct a $G$-space $\p$ with a natural $G$-projective-bundle on it, and shall show that the composite map
	\[[X;\p]_G\rightarrow\,{\rm Proj}_G(X)\rightarrow H^3_G(X;\Z)\]
is an isomorphism.

To prove the injectivity of Proj$_G(X)\rightarrow H^2_G(X;\T)\cong H^3_G(X;\Z)$ we consider the filtration
	\[{\rm Proj}_G(X)\supseteq\,{\rm Proj}^{(1)}\supseteq\,{\rm Proj}^{(0)},\]
where Proj$^{(1)}$ consists of the stable projective bundles which are trivial when the $G$-action is forgotten, i.e. those that can be described by cocycles
	\[\alpha:G\times X\rightarrow PU(\H)\]
such that
	\[\alpha(g_2,g,x)\alpha(g_1,x)=\alpha(g_2g_1,x),\]
and Proj$^{(0)}$ consists of those such that $\alpha$ lifts to
	\[\alpha :G\times X\rightarrow U(\H)\]
such that
	\be \alpha(g_2,g_1x)\alpha(g_1,x)=c(g_2,g_1,x)\alpha(g_2g_1,x)\ee
for some $c:G\times G\times X\rightarrow\T$.

We shall compare the filtration of Proj$_G(X)$ with the filtration
	\[H^2_G(X;\T)=H^{(2)}\supset H^{(1)}\supset H^{(0)}\]
defined by the spectral sequence.  By definition $H^{(1)}$ is the kernel of
	\[H^2_G(X;\T)\rightarrow E^{02}_1 = H^2(X;sh(\T))=\,{\rm Proj}(X),\]
and the composite
	\[{\rm Proj}_G(X)\rightarrow H^2_G(X;\T)\rightarrow \,{\rm Proj}(X)\]
is clearly the map which forgets the $G$-action.  Thus Proj$_G(X)/{\rm Proj}^{(1)}$ maps injectively to
	\[H^2_G(X;\T)/H^{(1)}\cong E^{02}_{\infty}\hookrightarrow E_1^{02}=\,{\rm Proj} (X).\]

Now let us consider the map
	\[{\rm Proj}^{(1)}\rightarrow H^{(1)}.\]
The subgroup $H^{(0)}$ is the kernel of $H^{(1)}\rightarrow E^{11}_2$, while $E^{11}_1=\,{\rm Pic}(G\times X)$.  We readily check that an element of Proj$^{(1)}$ defined by the cocycle
	\[\alpha :G\times X\rightarrow PU(\H)\]
maps to the element of Pic$(G\times X)$ which is the pull-back of the circle bundle $U(\H)\rightarrow PU(\H)$, and can conclude that $\alpha$ maps to zero in $E^{11}_2$ if and only if it defines an element of Proj$^{(0)}$.  Thus Proj$^{(1)}/{\rm Proj}^{(0)}$ injects into
	\[H^{(1)}/H^{(0)}=E^{11}_{\infty}=\,{\rm ker}:E^{11}_2\rightarrow E^{30}_2.\]
Finally, assigning to an element $\alpha$ of Proj$^{(0)}$ the class in
	\[E^{20}_2=H^2_{\rm c.c.}(G;{\rm Map}(X;\T))\]
of the cocycle $c$ occurring in (6.7), we see that if $[c]=0$ then the projective bundle comes from a $G$-Hilbert-bundle, which is necessarily trivial, as we have already explained.  So Proj$^{(0)}$ injects into $H^{(0)}=E_2^{20}$.
\bigskip

We now turn to the construction of the potential universal $G$-space $\mathcal P$ mentioned above.  We shall begin with a few general remarks about $G$-equivariant homotopy theory when $G$ is a compact group.

If $Y$ is a $G$-space we can consider the space $Y^H$ of $H$-fixed-points for any subgroup $H$ of $G$.  This is a space with an action of $W_H=N_H/H$, where $N_H$ is the normalizer of $H$ in $G$.  To give the space $Y^H$ clearly determines $[X;Y]_G$ when $X$ is a $G$-space of the form $X=(G/H)\times X_0$, where $G$ acts trivially on $X_0$; and to give $Y^H$ together with its $W_H$-action determines $[X;Y]_G$ whenever $X$ is isotypical of type $H$ (i.e. all isotropy groups in $X$ are conjugate to $H$), for then $[X;Y]_G$ is the homotopy classes of sections of a bundle on $X/G$ with fibre $Y^H$ associated to the principal $W_H$-bundle $X^H\rightarrow X/G$.

To give an element of Proj$_G(X)$ on an $H$-isotypical $G$-space $X$ is the same as to give a stable $N_H$-equivariant bundle on $X^H$. Because isomorphism classes of stable $H$-Hilbert-spaces correspond to elements of Ext$(H;\T)$, these bundles are classified by $W_H$-equivariant maps from $X^H$ to
	\[\p_H=\coprod_{\H\in\,{\rm Ext}(H;\T)}BPU(\H)^H,\]
where we represent an element of Ext$(H;\T)$ by the essentially unique Hilbert space $\H$ with a stable projective representation of $H$ inducing the extension.  The group $PU(\H)^H$ is disconnected, its group of components being Hom$(H;\T)$, but each connected component has the homotopy type of $B\T\cong\P^\infty_{\C}$.   As the classifying space functor $B$ commutes with taking $H$-invariants, the space $\p_H$, being a space of $H$-fixed-points, has a natural action of $W_H$. We shall give each group $PU(\H )$ the norm topology: there is then a natural projective bundle on $\p _H$ with fibres $\P(\H)$ which satisfies both conditions (i) and (ii) from the beginning of this section.

There is now a standard procedure --- unappealingly abstract --- for cobbling together a $G$-space $\p$  so that for each subgroup $H$ of $G$ we have $\p^H\simeq\p_H$.  We introduce the topological category $\mathcal O$ of $G$-orbits (i.e. transitive $G$-spaces) and $G$-maps.  Any $G$-space $Y$ gives a contravariant functor from $\mathcal O$ to spaces by
	\[S\mapsto \,{\rm Map}_G(S;Y).\]
If $S=G/H$, then Map$_G(G/H;Y)\cong Y^H$.  Conversely, suppose that $F$ is a contravariant functor from $\mathcal O$ to spaces.  Let $\mathcal O _F$ denote the topological category whose objects are triples $(S,s,y)$, where $S$ is an orbit, $s\in S$, and $y\in F(S)$.  A morphism $(S_0,s_0,y_0)\rightarrow (S_1,s_1,y_1)$ is a map $\theta :S_0\rightarrow S_1$ in $\mathcal O$ such that $\theta (s_0)=s_1$ and $\theta^*(y_1)=y_0$.  The group $G$ acts on the category $\mathcal O _F$ by
	\[g.(S,s,y)=(S,gs,y),\]
and so the ``realization"  $\vert\mathcal O _F\vert$ (in the sense of [S1]) is a $G$-space, and the fixed-point set $\vert\mathcal O _F\vert^H$ plainly contains $F(G/H)$. If each space $F(S)$ is an ANR then $\vert\mathcal O _F\vert$ is a $G$-ANR.

\setcounter{prop}{7}
\begin{prop} The inclusion $F(G/H)\rightarrow\vert\mathcal O_ F\vert^H$ is a homotopy-equivalence.\end{prop}

We shall omit the proof, which is quite elementary.  We apply it to the functor $F$ defined by $F(G/H)= \p _H$. There is no trouble in seeing that $\p=\vert\mathcal O _F\vert$ carries a tautological $G$-projective-bundle, so that we have a $G$-map
\setcounter{equation}{8}
	\be \p\rightarrow\,{\rm Map}(EG;BPU(\H)) \ee
into the space which represents the functor $X\mapsto H^3_G(X;\Z)$.  To see that (6.9) induces an isomorphism
	\[[X;\p]_G\rightarrow H^3_G(X;\Z)\]
it is enough (by the result of [JS]) to check the cases $X=(G/H)\times S^i$, when $S^i$ is an $i$-sphere; but this reduces to the isomorphism
	\[\pi_i(\p_H))\cong H^{3-i}(BH;\Z)\]
which we have already pointed out.

\section*{Appendix 1: The compact-open topology}
\renewcommand{\thesection}{A1}\setcounter{prop}{0}

The compact-open topology on the space Map$(X;Y)$ of continuous maps from a space $X$ to a metric space $Y$ is the topology of uniform convergence on all compact subsets of $X$.  (In fact there is no need for $Y$ to be metrizable, for the compact-open topology can also be defined as the coarsest topology for which the subsets
	\[F_{C,U}=\{f:X\rightarrow Y \ \ \mbox{such that} \ \ f(C)\subset U\}\]
are open whenever $C$ is compact in $X$ and $U$ open in $Y$.)  With this topology it is clear that a map $Z\rightarrow$ Map$(X;Y)$ is continuous if and only if the adjoint map $Z\times X\rightarrow Y$ is continuous on all subsets of the form $Z\times C$, where $C$ is compact in $X$.  If $Z$ and $X$ are metrizable this is simply saying that $Z\times X\rightarrow Y$ is continuous.

On the space Hom$({\mathcal H}_0;{\mathcal H}_1)$ of continuous linear maps between two Hilbert spaces the compact-open topology is only very slightly finer than the topology of pointwise convergence, which is called ``the strong operator topology" by functional analysts.  The Banach-Steinhaus theorem\footnote{Strictly, the Banach-Steinhaus theorem ([T] Thm 33.1, [B] chap.III \S 3,thm 2), which holds whenever ${\mathcal H}_0$ is Fr\'{e}chet and ${\mathcal H}_1$ is locally convex, asserts that a set of maps which is compact for the topology of pointwise convergence is equicontinuous.  But it is easy to see ([T] 32.5) that on equicontinuous  subsets the compact-open and pointwise topologies coincide.} tells us that exactly the same subsets are compact in these two topologies; and on compact subsets the topologies must of course coincide.  In particular, if $Z$ is a metrizable space the continuous maps $Z\rightarrow$ Hom$({\mathcal H}_0;{\mathcal H}_1)$ are the same for both topologies.
\bigskip

For a Hilbert space ${\mathcal H}$ the groups GL$({\mathcal H})$ and $U({\mathcal H})$ are subsets of End$({\mathcal H})$, but when we speak of the compact-open topology on these groups we mean their subspace topology not in End$({\mathcal H})$ but in End$({\mathcal H})\times$ End$({\mathcal H})$, in which they are embedded by $g\mapsto (g,g^{-1})$.  The reason is that on the subset $G$ of invertible elements of End$({\mathcal H})$ the map $G\rightarrow$ End$({\mathcal H})$ given by inversion is not continuous.  (For example, let $g_n$ be the diagonal transformation of the standard Hilbert space $l^2$ of sequences defined by
	\[\begin{array}{lclll} (g_n\xi)_k &=& \xi_k & \mbox{if} & k\neq n,\\ &=& n^{-1}\xi_n & \mbox{if} & k=n.\end{array}\]
then $g_n\xi\rightarrow\xi$ as $n\rightarrow\infty$ for every $\xi\in l^2$.  But if $\xi \in l^2$ is the vector with $\xi_k=k^{-1}$ then 
	\[\| g_n^{-1}\xi-\xi\|\rightarrow 1\]
as $n\rightarrow\infty$, so $g_n^{-1}\xi\not\rightarrow\xi$.)  Even when we define the compact-open topology so as to make inversion continuous, however, neither GL$({\mathcal H})$ and U$({\mathcal H})$ are quite topological groups, for the multiplication map is continuous only on compact subsets.  One can say that they are ``groups in the category of compactly generated spaces".  (See [St]. Functional analysts use the word {\it hypocontinuous} for bilinear maps which are continuous on compact subsets: the tensor product of distributions is a well-known example.)  In any case, for any metrizable space $Z$ the space of continuous maps into GL$({\mathcal H})$ or U$({\mathcal H})$ forms a group, and that is quite enough for our purposes.

We should also point out that the involution End$({\mathcal H})\rightarrow$ End$({\mathcal H})$ given by $A\mapsto A^*$ is not continuous for the compact-open topology.  For example let $A_n=e_0\otimes e^*_n$ be the operator of rank 1 in $l^2$ which takes $\xi =(\xi_k)$ to $A_n\xi=(\xi_n,0,0,0,...)$.  Clearly $A_n\rightarrow 0$ pointwise as $n\rightarrow\infty$.  But $A^*_n=e_n\otimes e^*_0$ takes the unit basis vector
	\[e_0=(1 \ 0 \ 0 \ 0 \ \ldots )\]
to the unit vector $e_n$, and so $A^*_ne_0\not\rightarrow 0$.
\bigskip

The most important positive result for our purposes is

\begin{prop} The group $U({\mathcal H})$ with the compact-open topology acts continuously by conjugation on the Banach space $\K(\H)$ of compact operators in $\H$, and also on the Hilbert space $\H^*\otimes\H$ of Hilbert-Schmidt operators. \end{prop}

\no {\bf Proof.}  (i) We must show that for each unitary operator $u_0$, each compact operator $k_0$, and each $\varepsilon >0$, we can find a compact subset $C$ of $\H$, and a $\delta >0$ such that if $\| k-k_0\|<\delta$ and $\| u(\xi)-u_0(\xi)\|<\delta$ for all $\xi\in C$ then
	\[\| uku^{-1}-u_0k_0u_0^{-1}\|<\varepsilon.\]
Now
	\begin{eqnarray*} \|uku^{-1}-u_0k_0u_0^{-1}\| &\leq & \|uku^{-1}-uk_0u^{-1}\|+\|uk_0u^{-1}-u_0k_0u^{-1}\|\\ && +\| u_0k_0u^{-1}-u_0k_0u_0^{-1}\|\\
&=&\|k-k_0\|+\|(u-u_0)k_0\|+\|k_0(u^*-u^*_0)\|\\ &=& \| k-k_0\|+\|(u-u_0)k_0\|+\|(u-u_0)k^*_0\|,\end{eqnarray*}
where in the last line we have used $\|A^*\|=\|A\|$.  Because $k_0$ and $k^*_0$ are both compact operators we can find a compact subset $C$ of $\H$ which contains $k_0\xi$ and $k^*_0\xi$ for all unit vectors $\xi$, and we get the desired inequality by taking $\delta =\varepsilon /3$.

(ii) If $k$ and $k_0$ are Hilbert-Schmidt operators, the preceding calculation remains true if the operator norms $\| \ \ \|$ are replaced by the Hilbert-Schmidt norm $\| \ \ \|_{HS}$, given by
	\[\|A\|^2_{HS}=\sum\|Ae_n\|^2,\]
where $\{e_n\}$ is an orthonormal basis of $\H$.  It is therefore enough to show that for any Hilbert-Schmidt $k_0$ we have
	\[\|(u-u_0)k_0\|_{HS}<\varepsilon\]
if $u-u_0$ is small in the compact-open topology.  But as $\|u-u_0\|<2$ we have
	\[\sum_{n>N}\|(u-u_0)k_0e_n\|^2\leq 4\sum_{n>N}\|k_0e_n\|^2,\]
which is $<\varepsilon /2$ for suitable $N$, and we can make
	\[\|(u-u_0)k_0e_n\|\]
small for all $n\leq N$.
\bigskip

That essentially completes our discussion of the compact-open topology, but we shall briefly mention a few other points.

Because a compact subset of End$(\H)$ is equicontinuous, it is bounded in the operator norm (even though the example of the sequence $\{e_0\otimes e^*_n\}$ above shows that the norm is not itself a continuous function).  This implies that $A\mapsto A^*A$ is continuous on compact sets, though $A\mapsto A^*$ is not.  Polynomial maps $A\mapsto p(A)$ are also continuous on compact sets, and hence --- as a continuous function on the spectrum can be uniformly approximated by polynomials --- so is the retraction map $A\mapsto (A^*A)^t$ used on the space of Fredholm operators in \S 3.

From the point of view of homotopy theory the one really bad feature of the compact-open topology is that the subspaces GL$(\H)$ and Fred$(\H)$ are neither open nor closed in the vector space End$(\H)$, and so are not ANRs.  In other words, if $X_0$ is a closed subspace of a space $X$ then a continuous map $X_0\rightarrow$ GL$(\H)$ need not be extendable to a neighbourhood of $X_0$ in $X$.

\section*{Appendix 2: Fredholm operators}
\renewcommand{\thesection}{A2}\setcounter{prop}{0}

\begin{prop} For a separable Hilbert space $\H$ the spaces ${\rm GL}(\H)$, ${\rm U}(\H)$, and ${\rm Fred}(\H)$ are contractible in the compact-open topology, by a homotopy
	\[h=\{h_t\}:X\times [0,1]\rightarrow X\] which is continuous on compact subsets. \end{prop}

\no {\bf Proof.}  A single map $h:$ End$(\H )\times [0,1]\rightarrow$ End$(\H)$ will deal with the three cases simultaneously: it will have the property that $h_t(g^{-1})=(h_t(g))^{-1}$, which is needed in view of the definition of the compact-open topology on GL$(\H)$ and U$(\H)$ which was explained in Appendix 1.

The essential point is that we can identify $\H$ with the standard Hilbert space $L^2([0,1])$ of complex-valued functions on the unit interval, and that then the projection operator $P_t$ which projects on to the first factor in
	\[L^2([0,1])=L^2([0,t])\oplus L^2([t,1])\]
depends continuously on $t\in [0,1]$ on the compact-open topology.  (For it is obviously continuous in the topology of pointwise convergence.)  Let us factorize $P_t$ as $i_tR_t$, where
	\[R_t:L^2([0,1])\rightarrow L^2([0,t])\]
is the restriction and $i_t$ is the inclusion of $L^2([0,t])$ in $L^2([0,1])$, and when $0<t\leq 1$ let us write
	\[Q_t:L^2([0,t])\rightarrow L^2([0,1])\]
for the isometric isomorphism given by
	\[(Q_tf)(x)=t^{1/2}f(tx).\]
Then we define $h_t:$ End$(\H)\rightarrow$ End$(\H)$ by
	\[h_t(A)=i_tQ_t^{-1}AQ_tR_t+(1-P_t)\]
when $t\in (0,1]$, and $h_0(a)=1$.  Because
	\[\|Q_tR_t\xi\|=\|P_t\xi\|\]
is continuous in $t$ and $\rightarrow 0$ as $t\rightarrow 0$, while
	\[\|i_tQ_t^{-1}A\|=\|A\|,\]
the homotopy $h_t$ from $h_1=$ (identity) to $h_0=$ (constant) is continuous as claimed, and it preserves the subsets GL$(\H)$, U$(\H)$, and Fred$(\H)$.

\begin{prop}  The space ${\rm Fred}'(\H)$ of Proposition 3.1 is a representing space for $K$-theory, i.e. for every compact space $X$ we have a natural bijection
	\[[X;{\rm Fred}'(\H)]\rightarrow K^0(X).\]\end{prop}

The proof, which follows closely the corresponding argument in the Appendix of $[A]$, will be presented as a sequence of lemmas in which we shall denote a map $X\rightarrow{\rm Fred}'(\H)$ by
	\[(A,B)=(\{A_x\},\{B_x\})_{x\in X},\]
where each $A_x$ is a Fredholm operator in $\H$ with parametrix $B_x$, and $A_xB_x-1$ and $B_xA_x-1$ depend continuously on $x$ in the {\it norm} topology.

\begin{lem} If $A_x$ is surjective (resp. injective) when $x=x_0$ then it is surjective (resp. injective) for all $x$ in a neighbourhood of $x_0$. \end{lem}

\no {\bf Proof.}  Suppose that $A_{x_0}$ is surjective.  Because the Fredholm operator $A_{x_0}B_{x_0}$ is of the form $1+$(compact) it has index 0, and so we can find a finite rank operator $F$ such that $A_{x_0}(B_{x_0}+F)$ is surjective, and hence an isomorphism.  As $A_x(B_x+F)$ depends continuously on $x$ in the norm topology, and invertible operators form an open set in the norm topology, we find that $A_x(B_x+F)$ is invertible for $x$ near $x_0$, and so $A_x$ is surjective there.  A similar argument applies when $A_{x_0}$ is injective.

\begin{lem}  Suppose that $A_x$ is surjective for all $x\in X$.  Then the spaces $E_x=\ {\rm ker}(A_x)$ form a finite dimensional vector bundle on $X$. \end{lem}

\no {\bf Proof.}  Given $x_0\in X$, let $\H_0=E^{\perp}_{x_0}$, and let $i_0:\H_0\rightarrow \H$ be the inclusion.  Then $A_x\circ i_0$ is bijective when $x=x_0$, and hence for all $x$ near $x_0$ by the preceding lemma.  Considering the map of short exact sequences
	\[ \begin{array}{rcrcl} \H_0 & \stackrel{i_0}{\longrightarrow} & \H & \longrightarrow &E_{x_0}\\ A_xi_0\downarrow && A_x\downarrow && \downarrow\\ \H & \longrightarrow & \H & \longrightarrow & 0\end{array}\]
we conclude that orthogonal projection defines an isomorphism $E_x\rightarrow E_{x_0}$ for all $x$ near $x_0$.

\begin{lem}  There is a subspace $\H_1$ of finite codimension in $\H$ such that $p\circ A_x$ is surjective for all $x\in X$, where $p$ is orthogonal projection $\H\rightarrow \H_1$. \end{lem}

\no {\bf Proof.}  By lemma A2.3 we can achieve this for $x$ in a neighbourhood of a chosen point of $X$.  But $X$ can be covered by a finite number of such neighbourhoods, and we can take the intersection of the corresponding subspaces $\H_1$.
\bigskip

\no {\bf Proof of Proposition A2.2.}  To each Fredholm family $(A,b)$ we can now associate the element
	\[\chi_{A,B}=[\{{\rm ker}(p\circ A_x)\}]-[X\times{\rm ker}(p)]\]
of $K^0(X)$, where $p$ is as in the preceding lemma.  The only choice made was of $\H_1$, but replacing $\H_1$ by a smaller subspace adds the same trivial bundle to both ker$(p\circ A)$ and $X\times$ ker$(p)$, so the $K$-theory class $\chi_{A,B}$, for a homotopy gives us an element of $K^0(X\times [0,1])\cong K^0(X)$.

Finally, we must show that if $\chi_{A,B}=0$ then $(A,B)$ is homotopic to a constant map.  But if $\chi_{A,B}=0$ we can assume (by making $\H_1$ smaller) that the bundle $\{{\rm ker}(p\circ A_x)\}$ is trivial, and isomorphic to $X\times$ ker$(p)$.  We can then add a finite rank family $\{F_x\}$ to $\{A_x\}$ so that $\tilde{A}_x=A_x+F_x$ is an isomorphism for all $x$; and $(\tilde{A},B)$ is still a map into Fred$'(\H)$, and is homotopic to $(A,B)$.  Because GL$(\H)$ is contractible in the compact-open topology, we can deform $(\tilde{A},B)$ to $(1,\tilde{A}^{-1}B)$, where $\tilde{A}^{-1}B$ is of the form $1+$ (compact), and then we can deform this family linearly to (1,1).

\section*{Appendix 3: Equivariant contractibility of the unitary group of Hilbert space in the norm topology}
\renewcommand{\thesection}{A3}\setcounter{prop}{0}

The results in this appendix are not, strictly speaking, needed in the paper, except to show that for a projective bundle with norm-continuous structure the two possible definitions of twisted equivariant $K$-theory coincide. We have included them partly for their intrinsic interest, and partly to correct a number of misstatements by the second author and others which have often been repeated in the literature.
\bigskip

Let $\H$ be a stable $G$-Hilbert-space, and $U(\H)$ the unitary group with the norm topology.  We have pointed out that the $G$-action on $\H$ does not induce a continuous action of $G$ on $U(\H)$.  The $G$-continuous elements $U_{G-{\rm cts}}(\H)=\{u\in U(\H):g\mapsto gug^{-1}$ is continuous$\}$ do, however, form a closed subgroup of $U(\H)$, in fact a sub-Banach-Lie-group. It is the intersection of $U(\H)$ with the closed linear subspace End$_{G-{\rm cts}}(\H)$ of End$(\H)$. To get a feeling for this subspace, notice that if $\H =L^2(G)$ then multiplication by an $L^{\infty}$ function $f$ on $G$ is a $G$-continuous operator if and only if $f$ is continuous. If $G$ is the circle group $\T$ then a $\T$-action on $\H$ defines a grading $\H = \oplus\H_k$, and any continuous linear map $A:\H\rightarrow \H$ can be represented by a block matrix $(A_{kl})$, where $A_{kl}:\H_l \rightarrow\H_k$. Roughly, $A$ is $G$-continuous if $||A_{kl}||\rightarrow 0$ sufficiently fast as $|k-l| \rightarrow \infty$

\begin{prop} The group $U_{G-{\rm cts}}(\H)$ is equivariantly contractible. \end{prop}

\begin{cor} The space Fred$_{G-{\rm cts}}(\H)$ of $G$-continuous Fredholm operators in $\H$, with the norm topology, is a representing space for $K^0_G$. \end{cor}

The corollary follows from the proposition by exactly the same argument used in the non-equivariant case in Appendix 2, and we shall say no more about it.

\bigskip
One can think of the results in the following way. Although $G$ does not act continuously on $U(\H)$ or Fred$(\H)$ it does make sense to say that a continuous map from a $G$-space $X$ to these spaces is $G$-equivariant. Then A3.1 says that any two $G$-maps $X \rightarrow U(\H)$ are homotopic, while A3.2 says that $K^0_G(X)$ is the set of homotopy classes of $G$-maps $X\rightarrow {\rm Fred}(\H)$. In this sense the misstatements referred to are innocuous.
\bigskip

\noindent {\bf Proof of (A3.1)}

Because $U=U_{G-{\rm cts}}(\H)$ is a $G$-ANR (see [JS]) it is enough to show that any $G$-map $f:X\rightarrow U$ from a compact $G$-space $X$ can be deformed to the constant map at the identity.  By a well-known ``Eilenberg swindle" argument it is enough to show that $f$ can be deformed into the subgroup of elements of the form
	\[\left(\begin{array}{cc} u & 0\\ 0 & 1\end{array}\right)\]
with respect to an orthogonal decomposition $\H=\H_1\oplus\H_2$ of $\H$ into stable $G$-Hilbert-spaces.  (For there is a canonical path from $u\oplus u^{-1}$ to the identity, and hence from
	\[u\oplus 1=u\oplus (1\oplus 1)\oplus (1\oplus 1)\oplus\ldots\]
to
	\[u\oplus (u^{-1}\oplus u)\oplus (u^{-1}\oplus u)\oplus \ldots\]
	\[=(u\oplus u^{-1})\oplus (u\oplus u^{-1})\oplus\ldots,\]
and hence to the identity.)

It is also enough if we perform the deformation in the larger group $GL=GL_{G-{\rm ctr}}(\H)$, for $GL$ can be equivariantly retracted to $U$ by the usual polar decomposition.
\bigskip

The essential step in Kuiper's proof is the

\begin{lem} For any $\varepsilon >0$ there is an orthogonal decomposition
	\[\H=\H_1\oplus\H_2\oplus\H_3\]
into stable $G$-Hilbert-spaces wuch that $f(x)(\H_1)$ is $\varepsilon$-orthogonal to $\H_3$ for every $x\in X$.  (We say that subspaces $P$ and $Q$ are $\varepsilon$-orthogonal if $\vert\langle p,q\rangle\vert<\varepsilon \| p\| \| q\|$ for all $p\in P$ and $q\in Q$.)\end{lem}

Granting the lemma, the proof of (A3.1) is as follows.  For each $x\in X$ we have an $\varepsilon$-orthogonal decomposition
\renewcommand{\theequation}{A3.\arabic{equation}}\setcounter{equation}{3}
	\be \H=f(x)\H_1\oplus\H_x\oplus\H_3,\ee
where $\H_x=\H\ominus (f(x)\H_1\oplus\H_3)$, and the projections on to each summand depend continuously on $x$ (in the norm topology).  Choose a fixed isomorphism $T:\H_1\rightarrow\H_2$.  Then the nearly unitary transformation $\varphi_x$ of $\H$ which, in terms of the decomposition (A3.4), takes
	\[f(x)\xi\oplus\eta\oplus T\zeta\]
to
	\[-f(x)\zeta \oplus\eta \oplus T\xi\]
belongs to $GL$, and is connected to the identity by the path obtained by conjugating the unitary rotation from
	\[\xi\oplus\eta\oplus\zeta \ \ \mbox{to} \ \ (-\zeta)\oplus\eta\oplus\xi\]
in $\H_1\oplus\H_2\oplus\H_1$ by $f(x)\oplus 1\oplus T$.  This path depends continuously on $x$.  The original map $f$ is therefore $G$-homotopic in $GL$ to $f_1$, where $f_1(x)=\varphi^{-1}_xf(x)$.  Now $f_1(x)\vert\H_1$ is simply the fixed map
	\[T:\H_1\rightarrow\H_3\subset \H,\]
so we can perform a rotation interchanging $\H_1$ and $\H_3$ to deform $f_1$ to a map $f_2$ such that $f_2(x)\vert\H_1$ is the identity for all $x\in X$.
\bigskip

\noindent {\bf Proof of Lemma (A3.3)}

Thinking of $f:X\rightarrow U_{G-{\rm ctr}}(\H)$ as a map into the Banach space End$(\H)$ we can find, because $X$ is compact, a map $\tilde{f}$ arbitrarily close to $f$ such that $\tilde{f}(X)$ is contained in a finite dimensional subspace $V$ of End$(\H)$.  In fact, because vectors $\xi\in\H$ with finite dimensional $G$-orbits are dense in $\H$ (cf. [CMS] p. 93), we can suppose $V$ is a $G$-subspace of $\H$, and, by averaging over $G$, that $\tilde{f}$ is a $G$-map, the image is automatically in $GL$.

Now suppose that we have found three orthogonal finite dimensional $G$-subspaces $P_1,P_2,P_3$ of $\H$ such that $P_1\cong P_3$ and $\alpha (P_1)\subset P_1\oplus P_2$ for all $\alpha\in V$.  Let $Q_1$ be an arbitrary irreducible $G$-subspace of $\H$ orthogonal to $P_1\oplus P_2\oplus P_3$.   We can clearly find two other finite dimensional subspaces $Q_2$ and $Q_3$, orthogonal both to each other and to $P_1\oplus P_2\oplus P_3\oplus Q_1$ so that $Q_1\cong Q_3$ and $\alpha (Q_1)\subset P_1\oplus Q_1\oplus P_2\oplus Q_3$ for all $\alpha\in V$.  Now define $P_i^{(1)}=P_i\oplus Q_i$ for $i=1,2,3$.  We have
	\[\alpha (P_1^{(1)})\subset P_1^{(1)}\oplus P_2^{(1)} \ \ \mbox{and} \ \ P_1^{(1)}\cong P_3^{(1)}.\]
Repeating the process we find increasing sequences of subspaces
	\[P_i\subset P_i^{(1)}\subset P_i^{(2)}\subset\ldots\]
such that $P_1^{(k)},P_2^{(k)},P_3^{(k)}$ are mutually orthogonal for all $k$, while
	\[\alpha (P_1^{(k)})\subset P_1^{(k)}\oplus P_2^{(k)} \ \ \mbox{and} \ \ P_1^{(k)}\cong P_2^{(k)}.\]
Finally we define $\H_1$ as the closure of the union of the subspaces $P_1^{(k)}$ for $k=1,2,...$, and $\H_3$ as the closure of the union of the $P_3^{(k)}$.  Then $\H_2$ is defined so that
	\[\H=\H_1\oplus\H_2\oplus\H_3.\]
It is obvious that we can make the choices so that all three subspaces $\H_i$ are stable.  We have now finished, for $\tilde{f}(x)(\H_1)$ is orthogonal to $\H_3$ for all $x\in X$, and so $f(x)(\H_1)$ is $\varepsilon$-orthogonal to $\H_3$ as $\| f(x)-\tilde{f}(x)\|<\varepsilon$.

\section*{References}

\begin{enumerate}
\item[[A]] Atiyah, M.F., K-theory.  Benjamin, New York, 1967.
\item[[AH]] Atiyah, M.F., and M.J.Hopkins, A variant of K-theory, $K_\Sigma$.  {\it In} Topology, Geometry, and Quantum Field Theory.  Ed. U. Tillmann, Cambridge Univ. Press, 2004.
\item[[AS]] Atiyah, M.F., and I.M.Singer, Index theory for skew-adjoint Fredholm operators.  Publ. Math. I.H.E.S. Paris, {\bf 37} (1969), 5--26.
\item[[B]] Bourbaki, N., Espaces vectoriels topologiques, chap. III--V.  Actualit\'{e}s Sci. et Ind. 1229, Hermann, paris 1954.
\item[[Bor]] Borel, A., Seminar on transformation groups. Ann. of Math. Stud. 46, Princeton Univ. Press, 1960.
\item[[BCMMS]] Bouwknegt, P., A.L.Carey, V.Mathai, M.K.Murray, D.Stevenson, Twisted $K$-theory and $K$-theory of bundle gerbes. Comm. Math. Phys. 228 (2002), 17--45.
\item[[CM]] Carey, A., and J.Mickelsson, The universal gerbe, Dixmier-Douady class, and gauge theory. Lett. Math. Phys., {\bf 59} (2002), 47--60.  
\item[[C]] Connes, A., Noncommutative geometry.  Academic Press, 1994.
\item[[CMS]] Carter, R., I.G.Macdonald, and G.B.Segal, Lectures on Lie groups and Lie algebras. London Math. Soc. Student Texts {\bf 32}, Cambridge Univ. Press, 1995.
\item[[D]] Dold, A., Partitions of unity in the theory of fibrations.  Ann. of Math., {\bf 78} (1963), 223--255.
\item[[DK]] Donovan, P., and M.Karoubi, Graded Braner groups and $K$-theory with local coefficients.  Publ. Math. I.H.E.S. Paris, {\bf 38} (1970), 5--25.
\item[[FHT]] Freed, D., Twisted $K$-theory and loop groups.  Proc. Int. Congress of Mathematicians Vol. III, Beijing, 2002, 419--430.
\item[[FGZ]] Frenkel, I.B., H.Garland, and G.J.Zuckerman, Semi-infinite cohomology and string theory.  Proc. Nat. Acad. Sci. U.S.A., {\bf 83} (1986), 8442--8446.
\item[[G]] Grothendieck, A., Le groupe de Brauer. S\'{e}minaire Bourbaki, Vol. 9, exp. 290, Soc. Math. France 1995.
\item[[HR]] Higson, N., and J.Roe, Analytic $K$-homology.  Oxford Univ. Press, 2000.
\item[[JS]] James, I.M., and G.B.Segal, On equivariant homotopy type.  Topology, {\bf 17} (1978), 267--272.
\item[[K]] Kuiper, N.H., The homotopy type of the unitary group of Hilbert space.  Topology, {\bf 3} (1965), 19--30.
\item[[Mi]] Mickelsson, J., Twisted $K$-theory invariants. Preprint AT/0401130.
\item[[Mo]] Moore, G., $K$-theory from a physical perspective. Topology, Geometry, and Quantum field theory. Ed. U.Tillmann. Cambridge Univ. Press, 2004.
\item[[P]] Palais, R.S., On the homotopy of certain groups of operators. Topology, {\bf 3} (1965), 271--279.
\item[[PS]] Pressley, A.N., and G.B.Segal, Loop groups.  Oxford Univ. Press, 1986.
\item[[R]] Rosenberg, J., Continuous-trace algebras from the bundle-theoretic point of view. J.Austral.Math.Soc., {\bf A47}(1989), 368-381.
\item[[S1]] Segal, G.B., Classifying spaces and spectral sequences.  Publ. Math. I.H.E.S., Paris, {\bf 34} (1968), 105--112\item[[S2]] Segal, G.B., Equivariant $K$-theory. Publ. Math. I.H.E.S., Paris, {\bf 34} (1968), 129-151. 
\item[[S3]] Segal, G.B., Fredholm complexes. Quarterly J. Math. Oxford, {\bf 21} (1970), 385--402.
\item[[S4]] Segal, G.B., Topological structures in string theory.  Phil. Trans. R. Soc. London, A{\bf 359} (2001), 1389--1398.
\item[[St]] Steenrod, N.E., A convenient category of topological spaces. Michigan Math. J.,{\bf 14}(1967),133-152.
\item[[T]] Tr\`{e}ves, F., Topological vector spaces, distributions, and kernels. Academic Press, 1967.
\item[[TXL]] Tu, J.-L., P. Xu, and C. Laurent-Gengoux, Twisted $K$-theory of differentiable stacks.  Ann. Sci. \'Ecole Norm. Sup. (4) {\bf 37} (2004), 841 -- 910.
\item[[W]] Wood, R., Banach algebras and Bott periodicity.  Topology, 4 (1965/6), 371--389.
\end{enumerate}

\end{document}